\documentclass{amsart}
\newtheorem{theorem}{Theorem}[section]

\newtheorem{corollary}[theorem]{Corollary}

\newtheorem{lemma}[theorem]{Lemma}

\renewcommand{\text}{\rm}

\begin{document}
\title{Some toric manifolds and a path integral}
\author{W. Rossmann}
\address{Department \ of Mathematics,
University of Ottawa\\
Ottawa, Cananda}
\email{rossmann@uottawa.ca}
\maketitle

\noindent{\small Preprint \ version of the paper published in Progress in
Math. 231, The Orbit Method in Geometry and Physics. In Honor of
A.A.Kirillov., \ Birkh\"auser, 2003, 395-419.}

\section{Introduction}

Toric manifolds and path integrals look like an odd couple at first sight,
but are in fact intimately related, a phenomenon having its origin in
statistical mechanics and quantum mechanics. It is not my intention to
discuss this topic any in depth -and I am ill equipped to do so in any event-
but only to suggest that this curiosity deserves study. It is in this spirit
that I have mixed some physical terminology into what is really a
mathematical exercise.

The toric manifolds in question are of a very special type. They were
invented by Bott and studied by Grossberg and Karshon \lbrack 1994\rbrack.
Their intriguing results stimulated my interest in these matters, especially
in the symplectic geometry of these manifolds. \ One thing that is new here
is the construction of 2-forms on these manifolds as curvature forms of
holomorphic line bundles with unitary metrics. These 2- forms are the objects
of interest and their realization here leads naturally to an explicit formula
from which everything else follows by simple calculations, things like
canonical coordinates and positivity criteria. The more important innovation,
however, seems to me the method, using a construction of these toric
manifolds as a special case of \ what I call '$\mathbb{P}^1$-chains', a
construct for which I have some hope. \ The toric example was worked out as a
test for the general case, to be discussed elsewhere. \ Here the methods are
completely elementary, relying on direct calculations presented in calculus
style.

A major reason for studying these manifolds comes from Kirillov's \lbrack
1968\rbrack\ orbit theory for unitary representations of Lie groups and
Kostant's \lbrack 1970\rbrack\ theory of quantization. \ The point of view
here is to look for the character of a unitary representation rather than for
the index of a virtual representation. This requires some vanishing theorems
for cohomology, which I deduce from a positivity result for the curvature by
Kodaira's theorem. The unitary structure on the representation space depends
on the curvature form and its positivity as well and the spectrum of the
representation is determined with the help of the K\"ahler potential and its
convexity. 

Kirillov's universal formula for the characters of representations
constructed by quantization has inspired much of what I have to say on this
topic here. In keeping with the physical origins of the subject, I have taken
the opportunity to explain in some detail two formulas for the character in
question -alias partition function- one along the lines of statistical
mechanics, the other in terms of a path integral. The first one is very
simple, but \ seems of some merit to me; the second one, the path integral
formula, seems of some interest as well, even though it may be commonplace in
a formal way. The point is simply to write out rather carefully how the path
integral can be defined in the situation at hand. \ The proof that it does
give the right answer is then immediate in any case, although the reason why
it does so -and not only here but also in many other situations- is far from
clear.

The paper consists then really of two parts, one dealing with the very
special toric varieties mentioned, the other with rather more general
questions, for which these toric varieties serve as vehicle of exposition. It
seemed to me that each part benefits from the other.

\section{$\boldsymbol{\mathbb{P}}^{\boldsymbol{1}}$ itself}

\noindent Write \text{G} for SL$(2,\mathbb{C})$ acting by holomorphic
automorphisms on $\mathbb{P}^1,$
\begin{equation*}
z(g\xi)=\frac{az(\xi)+c}{bz(\xi)+d},\quad ad-bc=1.
\end{equation*}
The points $\xi_0,\xi_\infty\,$ with $z=0,\infty$ are the fixed points of the
subgroup \text{T} of diagonal matrices and
$z=z(\xi)\in\mathbb{C}\cup\{\infty\}\,$is the affine coordinate on
$\mathbb{P}^1$ centered at $\xi_0.\,$ The weights by which \text{T} acts on
the cotangent spaces at these points $\xi_0,\xi_\infty$, denoted
$\alpha,-\alpha$, are twice the weights on the two coordinate lines
$\mathbb{C}e_{\alpha/2},\mathbb{C}e_{-\alpha/2}$ in $\mathbb{C}^2$ which
correspond to the points $\xi_0,\xi_\infty$ in $\mathbb{P}^1$. They are
interchanged by the reflection $s=s_\alpha$, transforming the affine
coordinate $z\,$centered at $\xi_0$ into the affine coordinate $z\circ
s=-1/z\,$centered at $\xi_\infty.\,$ Choose $\xi_0$ as base-point for
$\mathbb{P}^1$ and let \text{B} denote its isotropy group in \text{G}, thus
providing the realization of $\mathbb{P}^1$ as the coset space
\text{G}$/{\text{B}}$. The triple $({\text{G}},{\text{B}},{\text{T}})$
determines uniquely a basis $H_\alpha,E_\alpha,E_{-\alpha}$ of the Lie
algebra of \text{G} so that \text{T} is generated by $H_\alpha$,
\text{B}$\,$by $H_\alpha$ together with $E_\alpha,$ and 
\begin{equation*}
[H_\alpha,E_{\pm\alpha}]=\pm
2E_{\pm\alpha},\quad[E_\alpha,E_{-\alpha}]=H_\alpha.
\end{equation*}
Write \text{N} and \text{N}$^{-}$ for the subgroups generated by $E_\alpha$
and $E_{-\alpha}$, respectively. The affine coordinate $z=z(g{\text{B}})$ on
\text{G}$/{\text{B}}$ centered at $\xi_0$ is determined by
$\,v=e^{z(g)E_{-\alpha}}$ $\,$if $g=vcn$ according to
\text{G}$\,\doteq{\text{N}}^{-}{\text{T}}{\text{N}}$. \ The symbol
$\doteq\,$indicates that the decomposition is valid only for
$g{\text{B}}\neq\xi_\infty$. The compact real form \text{U}$\,:=$SU$(2)$ of
\text{G} has \text{i}$H_\alpha,E_\alpha\pm{\text{i}}E_{-\alpha}$ as basis for
its Lie algebra. $\,$

Homomorphisms \text{T}$\rightarrow$GL$(1,\mathbb{C})$ will be written
exponentially as $h\mapsto h^\lambda,\lambda$ being a linear functional on
the Lie algebra, referred to as a weight of \text{T}. A homomorphism of T
extends to \text{B}$\,={\text{T}}{\text{N}}\,$so that $b^\lambda=h^\lambda$
if $b=hn.\,\,$ $\mathbb{C}^2$ carries a symplectic form invariant under
\text{G} and a unitary metric invariant under \text{U}. Write $g\mapsto
g^\ast$ for the involution of \text{G} acting as $u^\ast=u^{-1}$ on \text{U}.
\ The decomposition \text{G}$={\text{U}}{\text{B}},\,g=ub,$ is valid
everywhere and unique if $b\in{\text{B}}={\text{T}}{\text{N}}$ is normalized
so that its component in \text{T} is real and positive.

The complex-valued holomorphic functions $f(g)\,$on \text{G} satisfying
$\,f(gb)=f(g)b^\lambda$ are the holomorphic sections of a line bundle
$\mathcal{L}_\lambda$ on $\mathbb{P}^1={\text{G}}/{\text{B}}$. In this
capacity $f(\xi)$ stands for the element in the line
$\mathcal{L}_\lambda(\xi)\,$at $\xi=g{\text{B}}$. This line bundle has a
\text{U} invariant unitary metric given by 
\begin{equation*}
|f(\xi)|^2:=|f(u)|^2\quad{\text{if}}\,\;\xi=u{\text{B}},\,u\in{\text{U}}.
\end{equation*}
One such function is $f(g):=\langle
e_{-\alpha\text{/2}},ge_{\alpha\text{/2}}\rangle^l$, the pointed brackets
denoting the \text{G}-invariant symplectic pairing on $\mathbb{C}^2$. The
curvature form $\sigma_\lambda$ of $\mathcal{L}_\lambda$, normalized so that
it represents the Chern class, is
\begin{equation*}
\sigma_\lambda=\frac{-\lambda(H_\alpha)}{2{\pi}{\text{i}}}\frac{d\overline{z}%
\wedge dz}{(1+\overline{z}z)^2}\,.
\end{equation*}
This $(1,1)$-form $\sigma_\lambda$ is \emph{negative} when
$\lambda(H_\alpha)$ is positive. These matters will be explained in a more
general setting when needed. 

In short, the triple $({\text{G}},{\text{B}},{\text{T}})$ determines the
basis $(H_\alpha,E_\alpha,E_{-\alpha})$, which allows the reconstruction of
$\mathbb{P}^1,$ complete with its affine coordinates $z\,,\,$the matrix
representation of its automorphism group on $\mathbb{C}^2$, its line bundles
$\mathcal{L}_\lambda$ with their unitary structures, etc. All of these
$\mathbb{P}^1$-paraphernalia are thus available for any triple of groups
known to be abstractly isomorphic with $({\text{G}},{\text{B}},{\text{T}})$.
Nevertheless, for the purpose of calculation in a more general setting it is
preferable to work directly with the triple
$({\text{G}},{\text{B}},{\text{T}})$. Three of its properties are assembled
here for reference. 

Let $v(g),c(g)$ be the components in \text{N}$^{-},{\,\text{T}}$ according to
the decomposition \text{G}$\,\doteq{\text{N}}^{-}{\text{T}}{\text{N}}$. Let
$u(g),b(g)$ be components in \text{U}, \text{B} according to
\text{G}$={\text{U}}{\text{B}}$, normalized so that $b(g)\in a(g){\text{N}}$
with $a(g)\in{\text{T}}$ real and positive.

\begin{lemma}
If $l:=\lambda(H_\alpha)=0,1,2,\cdots$ is a non-negative integer, then
$\,c_s^\lambda(g):=c(sg)^\lambda$ extends to a holomorphic function on all of
\text{G} and $\,c_s^\lambda(g)=0$ iff $g\in s{\text{B}}$.

\end{lemma}

\begin{proof}
In terms of the \text{G}-invariant pairing on $\mathbb{C}^2$, 
\begin{equation*}
c_s^{\alpha/2}(g)=\langle se_{_{\alpha/2}},ge_{\alpha/2}\rangle,\quad\langle
se_{\alpha/2},e_{\alpha/2}\rangle=1.
\end{equation*}
If $l:=\lambda(H_\alpha)$ then $\lambda=l(\alpha/2)$ and
\begin{equation*}
c_s^\lambda(g)=\langle se_{_{\alpha/2}},ge_{\alpha/2}\rangle^l.
\end{equation*}
The assertion follows.
\end{proof}

\begin{lemma}
\ If $u=u(g),\,v=v(g)$ for some $g$, then $c(u)=c^{-1}(b(v)\,)$.
\end{lemma}

\begin{proof}
\ If $u=u(g),\,v=v(g)$ then $u{\text{B}}=v{\text{B}}$ so $v=ub$ for some $b$
necessarily equal to $b(v).\,$Hence $u=vb^{-1}$ and
$c(u)=c(b^{-1})=c^{-1}(b)$.
\end{proof}

\noindent The third lemma is the essential one; its formulation may seem
awkward at this point.

\begin{lemma}
Let $e^{zE_{-\alpha}}\in{\text{N}}^{-},\,z=re^{{\text{i}}\phi},$ and
$a\in{\text{T }}$ real and positive. Then

\begin{equation*}
b(ae^{zE_{-\alpha}})\in\tilde{a}{{\text{N}}},\quad\tilde{a}:={a{%
\text{a}}}^{H_\alpha/2}(\tilde{r}),\quad\tilde{r}:=a^{-\alpha}r
\end{equation*}
where
\begin{equation*}
{\text{a}}(r):=1+r^2\text{ .}
\end{equation*}
and \text{a}$^{H_\alpha/2}:=e^{\frac{1}{2}(\log{\text{a}})H_\alpha}$.
\end{lemma}

\begin{proof}First take $a=1$ and set $v=e^{zE_{-\alpha}}$. Then
$ve_{\alpha/2}=e_{\alpha/2}+ze_{\alpha/2}$. Write $\,v=ub\,$and
$b\in\tilde{a}{\text{N}}$ with $\tilde{a}$ real and positive.$\,$Then
$ve_{\alpha/2}=\tilde{a}^{\alpha/2}ue_{\alpha/2}.\,$ Hence 
\begin{equation*}
\tilde{a}^{\alpha/2}ue_{\alpha/2}=e_{\alpha/2}+ze_{\alpha/2}
\end{equation*}
Taking square norms and using the normalization $\tilde{a}=\,$real, positive,
find 

\begin{equation*}
\tilde{a}^\alpha=1+r^2={\text{a}}(r).
\end{equation*}

\noindent Hence $\tilde{a}={{\text{a}}^{-H_\alpha/2}(z)}$, which is the
desired relation
$b(e^{zE_{-\alpha}})\in{\text{a}}^{H_\alpha/2}(z){\text{N}}\,$ for $a=1$. For
general $a$, 
\begin{equation*}
b(ae^{zE_{-\alpha}})=b(e^{a^{-\alpha}zE_{-\alpha}})a\in{\text{a}}^{H_%
\alpha/2}(a^{-\alpha}z)a{\text{N}}\,.
\end{equation*}
\end{proof}

\section{Definitions }

\noindent Let $(G_\ell,B_\ell,T_\ell),\,\cdots(G_1,B_1,T_1)$ be a sequence
obtained from standard $\mathbb{P}^1$-triples
$({\text{G}}_{\alpha_i},{\text{B}}_{\alpha_i},{\text{T}}_{\alpha_i})$ by
simply enlarging the 1-tori \text{T}$_i$ to bigger tori $T_i$ through central
extensions:
\begin{equation*}
G_i={\text{G}}_{\alpha_i}T_i,\quad
B_i={\text{B}}_{\alpha_i}T_i,\quad{\text{G}}_{\alpha_i}\cap
T_i={\text{T}}_{\alpha_i}={\text{B}}_{\alpha_i}\cap T_i\,.
\end{equation*}
All groups are required to be complex algebraic. The scheme is to denote by
roman letters objects belonging to $\mathbb{P}^1$ itself, like the triples
$({\text{G}}_\alpha,{\text{B}}_\alpha,{\text{T}}_\alpha).\,$ In addition, let
$T_i\rightarrow T_{i-1},h_i\mapsto h_i\gamma_{i\,i-1},$ be a sequence of
homomorphisms of these tori, written on the \emph{right}. The semidirect
decomposition
\text{B}$_{\alpha_i}={\text{T}}_{\alpha_i}{\text{N}}_{\alpha_i}$ extends
to$\,B_i=T_i{\text{N}}_{\alpha_i}$ and the homomorphisms extend to
\begin{equation*}
B_i={\text{N}}_{\alpha_i}T_i\rightarrow
T_i\overset{\gamma_{i\,i-1}}{\longrightarrow}T_{i-1}\subset B_{i-1},
\end{equation*}
producing the sequence
\begin{equation*}
G_\ell\overset{\subset}{\leftarrow}B_\ell\overset{\gamma_{\ell,\ell-1}}{%
\rightarrow}G_{\ell-1}\overset{\subset}{\leftarrow}\cdots\overset{\subset}{%
\leftarrow}B_2\overset{\gamma_{21}}{\rightarrow}G_1\,.
\end{equation*}
The $\mathbb{P}^1$-\emph{chain} associated to this sequence is the quotient
\begin{equation*}
\mathcal{X}:=G_\ell\times_{B_\ell}\times
G_{\ell-1}\times\cdots\times_{B_2}\times G_1/B_1
\end{equation*}

\noindent of $G_\ell\times G_{\ell-1}\times\cdots\times G_1$ by the (free)
right action of $B_\ell\times B_{\ell-1}\times\cdots\times B_1$ defined by
\begin{equation*}
(p_\ell,p_{\ell-1},\cdots,p_1)\mapsto(p_\ell
b_\ell,b_\ell^{-1}\gamma_{\ell\,\ell-1}p_{\ell-1},\cdots
p_2b_2,b_2^{-1}\gamma_{21}p_1).
\end{equation*}
The homomorphisms $\gamma$ are the \emph{connecting homomorphisms} of the
$\mathbb{P}^1$-chain. The construction can proceed one step at a time, from
right to left. Write \text{s}$\,:=$ $(G_\ell,B_\ell,T_\ell,$
$\gamma_{\ell\ell-1},\cdots,$ $\gamma_{21},G_1,B_1,T_1)$ for the whole
sequence of data which goes into the construction and
$\mathcal{X}=\mathcal{X}_{{\text{s}}}$ for the result. Write $\xi_0$ for the
element of $\mathcal{X}_{{\text{s}}}$ represented by $g_i=1$ for all $i$,
$\xi=g_\ell.g_{\ell-1}.\cdots.g_1.\xi_0$ for the element represented by
$(g_\ell,g_{\ell-1},\cdots,g_1)$ and $G_\ell.G_{\ell-1}.\cdots.G_1.\xi_0$ for
$\mathcal{X}_{{\text{s}}}$ itself$\,$. For brevity, write also $\xi$ as
$g_{\ell-1}.g_\ell\cdots$ $\,$ and $\mathcal{X}_{{\text{s}}}$ as
$G_\ell.G_{\ell-1}\cdots$. Equality in the quotient
$\mathcal{X}_{{\text{s}}}$ means
\begin{equation*}
g_\ell
b_\ell.g_{\ell-1}\cdots=g_\ell.b_\ell\gamma_{\ell\,\ell-1}g_{\ell-1}\cdots.
\end{equation*}
Omit $\gamma_{\ell\,\ell-1}$ and write $b_\ell g_{\ell-1}$ for
$b_\ell\gamma_{\ell\,\ell-1}g_{\ell-1}$ when the meaning is clear from the
indices. Let
$G_{{\text{s}}}:=\prod_{{\text{s}}}g_i,B_{{\text{s}}}:=\prod_{{\text{s}}}B_i$
denote the direct products and
$\mathcal{X}_{{\text{s}}}:=G_{{\text{s}}}/_\gamma B_{{\text{s}}},$ the
quotient by the right action of $B_{{\text{s}}}$ on $G_{{\text{s}}}\,$defined
above. The product torus $T_{{\text{s}}}:=\prod_{{\text{s}}}T_i$ acts on
$\mathcal{X}_{{\text{s}}}$ on the left, denoted $\xi\mapsto h\xi$ and given
by 
\begin{equation*}
g_\ell.g_{\ell-1}.\cdots\mapsto h_\ell g_\ell.h_{\ell-1}g_{\ell-1}.\cdots.
\end{equation*}
\noindent$T_{{\text{s}}}$ has an open dense orbit on
$\mathcal{X}_{{\text{s}}}$, so that
$(\mathcal{X}_{{\text{s}}},T_{{\text{s}}})$ is a smooth toric variety in the
usual sense.

Associated to the sequence $(G_i,B_i,T_i,\gamma_{ii-1})\,$is the sequence of
roots $\alpha_i\,$and a system of integers $c_{ji},j>i,$ defined as follows.
Let $\,h_j\rightarrow h_j\gamma_{ji}$ be the composite of the connecting
homomorphisms from $T_j$ to $T_i$ and let
$\gamma_{\,ji}\alpha_i\leftarrow\alpha_i$ be the pull-back of the root to
$T_j$ from $T_i$ . The value of $\gamma_{\,ji}\alpha_i$ on the canonical
generator $H_{\alpha_j}$ of the subtorus \text{T}$_{\alpha_j}$ of $T_j$ is an
integer denoted $c_{ji}:=\langle H_{\alpha_j},\gamma_{ji}\alpha_i\rangle$.
The homomorphism $\gamma_{ji}$ will be omitted if the meaning of $\langle
H_{\alpha_j},\alpha_i\rangle$ is clear from the indices. 

\medskip Some comments. The example which motivated the definition is the
case when the $(G_i,B_i,T_i)$ are triples naturally associated to a sequence
of roots $\alpha_i$ for a triple $(G,B,T)$ consisting of a semisimple group
$G$, a Borel subgroup $B$ and a Cartan subgroup $T$. In this case $T_i=T$ for
all $i$ and one can take the identity map $T\rightarrow T$ to construct the
connecting homomorphisms $B_i\rightarrow B_{i-1}$. If one replaces the
triples $(G_i,B_i,T)$ by triples $(P_i,B,T)$ where the $P_i$ are minimal
(proper) parabolic subgroups of $G$ containing $B$ and takes the identity
$B\rightarrow B$ for connecting homomorphisms one obtains the well-known
Bott-Samelson manifolds. Other connecting homomorphisms, which are generally
an essential ingredient, are of interest as well, e.g. those induced by
conjugations in $G$ for the triples $(P_i,B,T)$. \ Manifolds isomorphic to
the toric $\mathbb{P}^1$-chains defined above were first defined by Bott
using a different procedure and studied by Grossberg and Karshon \lbrack
1994\rbrack, as mentioned in the introduction. Although it would not be worth
while introducing $\mathbb{P}^1$-chains constructed by means of connecting
homomorphisms just for the case discussed here, there is a generalization for
which it is: each triple $(G_i,B_i,T_i)$ can be taken to consist of an
arbitrary algebraic groups $G_i$ whose quotient $G_i/B_i$ by its maximal
solvable subgroup is isomorphic with $\mathbb{P}^1$. \ An extension of
SL$(2,\mathbb{C})$ by a Heisenberg-Weyl group is an example, beyond those
already mentioned. 

\ 

\section{The curvature formula}

\noindent Let $\,\mathcal{X}_{{\text{s}}}$ be the $\mathbb{P}^1$-chain built
from a sequence
\begin{equation*}
{\text{s}}:=(G_\ell,B_\ell,T_\ell\,,\gamma_{\ell\,\ell-1},\,\cdots\,%
\gamma_{21},G_1,B_1,T_1\,).
\end{equation*}
\textbf{Line bundles}. Let $c(g_i)$ be the component in $T_i$ according to
the decomposition $G_i\doteq{\text{N}}^{-}_{\alpha_i}T_iN_{_i}$, defined for
$g_i$ in this open set. \ A weight $\lambda$ of the product torus
$T_{{\text{s}}}$ is a sequence $\lambda=(\lambda_i)$. Define a function
$c_{{\text{s}}}^\lambda\,$on the open set $\prod_{{\text{s}}}B_is_iB_i\,$in
$G_{{\text{s}}}$ by the formula
\begin{equation*}
c_{{\text{s}}}^\lambda(g_{_\ell},\cdots,g_1):=\prod_{{\text{s}}}c_{s_i}^{%
\lambda_i}(g_i),\quad c_{s_i}(g_i):=c(s_ig_i).
\end{equation*}
This function is defined on the domain where all the factors are defined,
i.e. $g_i\in{s_iN^{-}_i}B_i=B_is_iB_i.\,$ We record a transformation property
of the function $c_{{\text{s}}}^\lambda$ under the right action of
$B_{{\text{s}}}$ on
$G_{{\text{s}}}.\,\,\,$Let$\,\gamma_{i\,i-1}\lambda_{i-1}$ denote the
pull-back of a weight $\lambda_{i-1}$ from $B_{i-1}$ to$B_i$ and set
\begin{equation*}
\gamma\lambda:=(\gamma_{\ell\,\ell-1}\lambda_{\ell-1},\cdots\gamma_{21}%
\lambda_1,0),\quad\,{\text{s}}\lambda:=(s_\ell\lambda_\ell,\cdots,s_1%
\lambda_1).\,
\end{equation*}
Also set $w_i:=s_i\cdots s_1$ and
\text{w}$\lambda:=(w_\ell\lambda_\ell,\cdots,w_1\lambda_1)\,$.

\begin{lemma}

As function of $g\in G_{{\text{s}}}$, $c_{{\text{s}}}^\lambda(g)$ transforms
under the right action of $b\,$according to the
rule$\,\,c_{{\text{s}}}^\lambda(g\cdot
b)=c_{{\text{s}}}^\lambda(g)b^{\tilde{\lambda}}\,$where 
\begin{equation*}
\tilde{\lambda}:=(1-{\text{s}}\gamma)\lambda,\,\text{
}\lambda=(1+{\text{w}}\,\gamma)\tilde{\lambda}\text{ .}
\end{equation*}

\end{lemma}

\begin{proof}
Omit the connecting homomorphisms from the notation to calculate
\begin{equation*}
\begin{aligned}[t]
c_{{\text{s}}}^\lambda(g\cdot b)&=c_{{\text{s}}}^\lambda(g_\ell
b_{\ell-1},b_\ell^{-1}g_{\ell-1}b_{\ell-1},\cdots b_2^{-1}g_1b_1)\\
&=c^{\lambda_\ell}(s_\ell g_\ell
b_\ell)\,c^{\lambda_{\ell-1}}(s_{\ell-1}b_\ell^{-1}g_{\ell-1}b_{\ell-1})%
\cdots\,c^{\lambda_1}(s_1b_2^{-1}g_1b_1)\\
&=c^{\lambda_\ell}(s_\ell
g_\ell)b_\ell^{\lambda_{\ell-1}}b_\ell^{-s_{\ell-1}\lambda_{\ell-1}}c^{%
\lambda_{\ell-1}}(g_{\ell-1})b_{\ell-1}^{\lambda_{\ell-1}}\cdots\,b_2^{-s_1%
\lambda_1}c^{\lambda_1}(g_1)b_1^{\lambda_1}\\
&=c_{{\text{s}}}^\lambda(g)b^{\tilde{\lambda}}
\end{aligned}
\end{equation*}
One checks that $\tilde{\lambda}:=(1-{\text{s}}\gamma)\lambda$ is equivalent
to $\lambda=(1+{\text{w}}\,\gamma)\tilde{\lambda}$ .\quad
\end{proof}

The locally defined holomorphic functions $f(g)$ on $G_{{\text{s}}}$
satisfying $f(g\cdot b)=f(g)b^{\tilde{\lambda}}$ where defined are the local
sections of a holomorphic line bundle $\mathcal{L}_\lambda$ on
$\mathcal{X}_{{\text{s}}}$. Write $f(\xi)$ for the value of the section $f$
in the line $\mathcal{L}_\lambda(\xi)\approx\mathbb{C}$ attached
$\xi\in\mathcal{X}_{{\text{s}}}$. This line $\mathcal{L}_\lambda(\xi)$
carries a unitary metric defined by $|f(\xi)|^2:=|f(u)|$ if $\xi=u.\xi_0$
with $u\in U_{{\text{s}}}:=\prod_{{\text{s}}}U_i$ $,$ $U_i$ being the compact
form of $G_i$. To make use of fiber bundle technology is useful to note this
reformulation of the lemma.

\begin{corollary}
$\mathcal{L}_\lambda=G_{{\text{s}}}\times_{B_{{\text{s}}}}\mathbb{C}_{\tilde{%
\lambda}}$ is the holomorphic line bundle on
$\mathcal{X}_{{\text{s}}}=G_{{\text{s}}}/_\gamma\,B_{{\text{s}}}$ associated
to the 1-dimensional representation$\,\tilde{\lambda}$
$=(1-{\text{s}\gamma)\lambda}\,$of the structure group $B_{{\text{s}}}$ of
the principal bundle $G_{{\text{s}}}\rightarrow\mathcal{X}_{{\text{s}}}$.
$\,$ On the underlying real manifold
$\mathcal{X}_{{\text{s}}}=U_{{\text{s}}}/_\gamma
T_{{\text{s}}{\text{U}}}\,(T_{{\text{s}}{\text{U}}}$ $:=T_{{\text{s}}}\cap
U_{{\text{s}}})$ the line bundle is
$\mathcal{L}_\lambda=U_{{\text{s}}}\times_{T_{{\text{s}}{\text{U}}}}%
\mathbb{C}_{\tilde{\lambda}}$ and its unitary metric is associated to the
$T_{{\text{s}}{\text{U}}}$ invariant unitary metric on
$\mathbb{C}_{\tilde{\lambda}}$.\hfill{\qedsymbol}
\end{corollary}

Any holomorphic line bundle $\mathcal{L}$ equipped with a unitary structure
has a $(1,1)$-form $\sigma$ as curvature. Normalized so that it represents
the Chern class, this form is given by
\begin{equation*}
\sigma=-\frac{1}{2{\pi}{\text{i}}}\bar{\partial}\partial\log|f|^2
\end{equation*}
for any local holomorphic section $f$, the formula being valid where $f\neq
0$ \lbrack Giffiths and Harris, 1978, p.142\rbrack. It will be computed
explicitly for the case at hand.

\medskip

\noindent\textbf{Affine coordinates}. The equation
\begin{equation*}
\xi:=e^{z_\ell
E_{\alpha_\ell}}s_\ell.e^{z_{\ell-1}E_{\alpha_{\ell-1}}}.s_{\ell-1}\,\cdots\,%
e^{z_1E_{\alpha_1}}s_1.\xi_0.
\end{equation*}
defines holomorphic functions $z_\ell,z_{\ell-1}$,$\cdots,z_1\,$on the image
of the dense open set $\prod_{{\text{s}}}B_is_iB_i$ under the quotient map
$G_{{\text{s}}}\rightarrow\mathcal{X}_{{\text{s}}}$. \ Call
$(z_\ell,\cdots,z_1)$ the \emph{affine coordinates} on
$\mathcal{X}_{{\text{s}}}$ centered at
$\xi_\infty:=s_\ell.s_{\ell-1}\,\cdots\,s_1.\xi_0\,.\,$These functions
$z_i=z_i(\xi)$ are well-defined because
$G_{{\text{s}}}\rightarrow\mathcal{X}_{{\text{s}}}$ maps
$\prod_{{\text{s}}}{\text{N}}_is_i\,$isomorphically onto its image, but their
definition requires that $\xi$ be written in the form indicated.

\begin{lemma}
Under the action of $\,$an element $h=(h_i)$ of $T_{{\text{s}}}$ the affine
coordinates $(z_i)$ transform according to the rule 
\begin{equation*}
z_i\circ h=h^{\varpi_i}z_i,\quad\varpi_i:=(w_{\ell
i}\alpha_i,\cdots,w_{i+1\,i}\alpha_i).
\end{equation*}
$w_{ji}:=s_j\gamma_{j\,j-1}s_{j-1}\cdots\gamma_{i+1i}s_i\,$ is the
homomorphism $T_j\rightarrow T_i$ composed of the reflections $s$ and the
connecting homomorphisms $\gamma$.

\end{lemma}

\begin{proof}
Let $\xi:=\exp(z_\ell
E_{\alpha_\ell})s_\ell.\exp(z_{\ell-1}E_{\alpha_{\ell-1}}).s_{\ell-1}\,\cdots%
\,$. Omit the connecting homomorphisms from the notation to calculate 
\begin{equation*}
\begin{aligned}[t]
h\xi&:=h_\ell\exp(z_\ell
E_{\alpha_\ell})s_\ell.h_{\ell-1}\exp(z_{\ell-1}E_{\alpha_{\ell-1}}).s_{%
\ell-1}\,\cdots\,\\
&=\exp(h_\ell^{\alpha_\ell}z_\ell E_{\alpha_\ell})h_\ell
s_\ell.h_{\ell-1}\exp(z_{\ell-1}E_{\alpha_{\ell-1}}).s_{\ell-1}\,\cdots\,\\
&=\exp(h_\ell^{\alpha_\ell}z_\ell
E_{\alpha_\ell})s_\ell.h_\ell^{s_\ell}h_{\ell-1}\exp(z_{\ell-1}E_{\alpha_{%
\ell-1}}).s_{\ell-1}\,\cdots\,\\
&=\exp(h_\ell^{\alpha_\ell}z_\ell
E_{\alpha_\ell})s_\ell.\exp(h_\ell^{s_\ell\alpha_{\ell-1}}h_{\ell-1}^{%
\alpha_{\ell-1}}z_{\ell-1}E_{\alpha_{\ell-1}})h_\ell^{s_\ell}h_{\ell-1}.s_{%
\ell-1}\,\cdots
\end{aligned}
\end{equation*}
etc. This gives 
\begin{equation*}
z_i\circ h=(h_\ell^{w_{\ell
i}\alpha_i}h_{\ell-1}^{w_{\,\ell-1i}\alpha_i}\cdots
h_{i+1}^{w_{\,i+1\,i}\alpha_i}h_i^{\alpha_i})z_i,
\end{equation*}
and hence the assertion.
\end{proof}

The points in the coordinate domain with all coordinates non-zero form a
single orbit of $T_{{\text{s}}}$, namely its open dense orbit on the toric
manifold $\mathcal{X}_{{\text{s}}}$. 

\medskip

\noindent\textbf{Curvature}. Let $\lambda=(\lambda_i)$ be a weight of
$T_{{\text{s}}}$ and $\mathcal{L}_\lambda$ the associated holomorphic line
bundle on $\mathcal{X}_{{\text{s}}}$ with its unitary metric. The following
formula gives its curvature.

\begin{theorem}
The curvature form $\sigma_\lambda$ of $\mathcal{L}_\lambda$ is given by
\begin{equation*}
\begin{aligned}[t]
&\sigma_\lambda=+\frac{1}{2{\pi}{\text{i}}}\sum_{ij}\frac{\partial^2K_%
\lambda}{\partial\bar{z}_i\partial z_j}\,d\bar{z}_i\wedge dz_j\,,\quad
K_\lambda:=\log a_\lambda,\\
&a_\lambda:={\text{a}}(\tilde{r}_1)^{l_1}\cdots{\text{a}}(\tilde{r}_\ell)^{l_%
\ell},\quad
l_i:=\lambda_i(H_{\alpha_i}),\quad{\text{a}}(\tilde{r}):=1+\tilde{r}^2.
\end{aligned}
\end{equation*}
$\tilde{r}_\ell,\cdots,\tilde{r}_1\,$ are the functions of
$r_\ell,\cdots,r_1$ and v.v. defined inductively by 
\begin{equation*}
r_\ell=\tilde{r}_\ell,\quad
r_i=\tilde{r}_i{\prod_{j:j>i}{\text{a}}^{c_{ji}/2}}(\tilde{r}_j).\,
\end{equation*}
The exponents $c_{ji}$ are the integers $c_{ji}=\langle
H_{\alpha_j},\alpha_i\rangle$. 
\end{theorem}

\noindent The variables $\tilde{r}_i$ are defined where the
$z_i=r_ie^{\sqrt{-1}\phi_i}$ are and $(\tilde{r}_i,\,\phi_i)\,$ may be used
as real coordinates on $\mathcal{X}_{{\text{s}}}$ instead of $(r_i,\phi_i)$.
Their recursive definition proceeds in the direction opposite to the
construction of the $\mathbb{P}^1$- chain. \ Combine the
$(\tilde{r}_i,\,\phi_i)$ into complex valued variables
$\tilde{z}_i:=\tilde{r}_ie^{\sqrt{-1}\tilde{\phi}_i},\,$considered as
functions of the $z_i=r_ie^{\sqrt{-1}\phi_i}$ by the above equations
supplemented by $\phi_i=\tilde{\phi}_i$, even though these functions
$\tilde{z}_i$ are not holomorphic. The angle variables
$\phi_i=\tilde{\phi}_i$ play no role in much of what is to follow and will
then be suppressed. (The hoary notion of 'variables', which leaves it to
reader and context to decide what is a function of what and in what way, is
most convenient here and even more so later.) 

\begin{proof}To simplify the notation take a $\mathbb{P}^1$-chain of length
three, and take $(\alpha,\beta,\gamma)$ as indices. Use the shorthand
notation like $\xi=g_\alpha.g_\beta.g_\gamma.\xi_0$ etc, and omit the
connecting homomorphisms. By definition
\begin{equation*}
|c_{{\text{s}}}^\lambda(g_\alpha.g_\beta.g_\gamma.\xi_0)|^2=|c^{\lambda_%
\alpha}(u_\alpha)|^2|c^{\lambda_\beta}(u_\beta)|^2|c^{\lambda_\gamma}(u_%
\gamma)|^2,
\end{equation*}
if $g_\alpha.g_\beta.g_\gamma.\xi_0=u_\alpha.u_\beta.u_\gamma.\xi_0.\,\,$In
terms of the components $u(g)\in{\text{U}},b(g)\in{\text{B}}\,$according to
\text{G}$\,={\text{U}}{\text{B}}$, 
\begin{equation*}
\begin{aligned}[t]
&s_\alpha g_\alpha.s_\beta g_\beta.s_\gamma g_\gamma.\xi_0=\\
&=u_\alpha b_\alpha.s_\beta g_\beta.s_\gamma
g_\gamma.\xi_0\quad[u_\alpha:=u(s_\alpha g_\alpha),b_\alpha:=b(s_\alpha
g_\alpha)]\\
&=u_\alpha.u_\beta b_\beta.s_\gamma g_\gamma.\xi_0\quad[u_\beta:=u(b_\alpha
s_\beta g_\beta),b_\beta:=b(b_\alpha s_\beta g_\beta)]\\
&=u_\alpha.u_\beta.u_\gamma b_\gamma.\xi_0\quad[u_\gamma:=u(b_\beta s_\gamma
g_\gamma),b_\gamma:=b(b_\beta s_\gamma g_\gamma)]\\
&=u_\alpha.u_\beta.u_\gamma.\xi_0
\end{aligned}
\end{equation*}
In terms of the components $v(g)\in{\text{N}}^{-}\,\,$according to
\text{G}$\,\doteq{\text{N}}^{-}{\text{T}}{\text{N}}$ (Lemma 2.2,
$b^\lambda:=c(b)^\lambda$), \ 
\begin{equation*}
\begin{aligned}[t]
|c^{\lambda_\alpha}\left(u_\alpha\right)|^2&=b^{-2\lambda_\alpha}(\tilde{v}_%
\alpha)\quad[\tilde{v}_\alpha:=v(s_\alpha g_\alpha)]\\
|c^{\lambda_\beta}\left(u_\beta\right)|^2&=b^{-2\lambda_\beta}(\tilde{v}_%
\beta)\quad[\tilde{v}_\beta:=v(b_\alpha s_\beta g_\beta)]\\
|c^{\lambda_\gamma}\left(u_\gamma\right)|^2&=b^{-2\lambda_\gamma}(\tilde{v}_%
\gamma)\quad[\tilde{v}_\gamma:=v(b_\beta s_\gamma g_\gamma)].
\end{aligned}
\end{equation*}
Define $\tilde{z}_\alpha,\tilde{z}_\beta,\tilde{z}_\gamma$ by
\begin{equation*}
e^{\tilde{z}_\alpha E_{-\alpha}}:=\tilde{v}_\alpha,e^{\tilde{z}_\beta
E_{-\beta}}:=\tilde{v}_\beta,e^{\tilde{z}_\gamma
E_{-\gamma}}:=\tilde{v}_\gamma.
\end{equation*}
so that
\begin{equation*}
|c_{{\text{s}}}^\lambda(g_\alpha.g_\beta.g_\gamma.\xi_0)|^2=b^{-2\lambda_%
\alpha}(e^{\tilde{z}_\alpha
E_{-\alpha}})\,b^{-2\lambda_\beta}(e^{\tilde{z}_\beta
E_{-\beta}})\,b^{-2\lambda_\gamma}(e^{\tilde{z}_\gamma E_{-\gamma}}).
\end{equation*}
By Lemma 2.3 with $a=1$ this becomes 
\begin{equation*}
|c_{{\text{s}}}^\lambda(g_\alpha.g_\beta.g_\gamma.\xi_0)|^2={\text{a}}^{-l_%
\alpha}(\tilde{z}_\alpha){\text{a}}^{-l_\beta}(\tilde{z}_\alpha){%
\text{a}}^{-l_\gamma}(\tilde{z}_\gamma).
\end{equation*}
It remains to calculate the
$\tilde{z}_\alpha,\tilde{z}_\beta,\tilde{z}_\gamma$ as functions of the
coordinates $z_\alpha,z_\beta,z_\gamma$ of 
\begin{equation*}
\xi:=e^{z_\alpha E_\alpha}s_\alpha.e^{z_\beta E_\beta}s_\alpha.e^{z_\gamma
E_\gamma}s_\gamma\xi_0\,.
\end{equation*}
Thus take
\begin{equation*}
g_\alpha:=e^{z_\alpha E_\alpha}s_\alpha,\,g_\beta:=e^{z_\beta
E_\beta}s_\alpha,\,g_\gamma:=e^{z_\gamma E_\gamma}
\end{equation*}
or equivalently
\begin{equation*}
s_\alpha g_\alpha=e^{-z_\alpha E_{-\alpha}},s_\alpha g_\beta=s_\alpha
e^{-z_\beta E_{-\beta}},s_\gamma g_\gamma=e^{-z_\gamma E_{-\gamma}}.
\end{equation*}
The problem is this. Given $b_\alpha\in B_\alpha$ find $b_\beta\in B_\beta$
from $b_\beta=b(b_\alpha e^{-z_\beta E_{-\beta}})$ in order to determine
$\tilde{v}_\gamma=e^{\tilde{z}_\gamma E_{-\gamma}}$
$\in{\text{N}}^{-}_\gamma$ from $\tilde{v}_\gamma=v(b_\beta e^{-z_\gamma
E_{-\gamma}})$.\smallskip

\noindent At the first step of the process determine $b_\alpha$ $\,$from
$b_\alpha=b(g_\alpha)$. For $g_\alpha=e^{-z_\alpha E_{-\alpha}}$ the equation
becomes
\begin{equation*}
b_\alpha=b(e^{-z_\alpha E_{-\alpha}}).
\end{equation*}
By Lemma 2.3 (which is not affected by the minus sign in the exponent)
\begin{equation*}
b_\alpha=b(e^{-z_\alpha
E_{-\alpha}})\in\tilde{a}_\alpha{{\text{N}}_\beta},\,\tilde{a}_\alpha={%
\text{a}}^{H_\alpha/2}(\tilde{z}_\alpha)\quad\tilde{z}_\alpha:=z_\alpha.
\end{equation*}
At the second step, determine $b_\beta\,$ from the equation
\begin{equation*}
\,b_\beta=b(b_\alpha e^{-z_\beta E_{-\beta}}),\quad
b_\alpha\in\tilde{a}_\alpha{\text{N}}_\alpha.
\end{equation*}
$b_\alpha e^{-z_\beta E_{-\beta}}$ is the action of $b_\alpha$ on $B_\beta$
through the connecting homomorphism, which is trivial on \text{N}$_\alpha$.
Hence $b_\alpha$ may be replaced by $\tilde{a}_\alpha$ and same lemma
applies: 
\begin{equation*}
\begin{aligned}[t]
b_\beta&=b(\tilde{a}_\alpha e^{-z_\beta
E_{-\beta}})\in\tilde{a}_\beta{\text{N}}_\beta,\\
\tilde{a}_\beta&:=\tilde{a}_\alpha{\text{a}}^{H_\beta}\,(\tilde{z}_\beta\,)={%
\text{a}}^{H_\alpha}(\tilde{z}_\alpha){\text{a}}^{H_\beta}\,(\tilde{z}_\beta\,%
),\\
\tilde{z}_\beta&:=(\tilde{a}_\alpha)^{-\beta}z_\beta={\text{a}}^{-\beta(H_%
\alpha)/2}(\tilde{z}_\alpha)z_\beta.
\end{aligned}
\end{equation*}
At the third step, 
\begin{equation*}
\,b_\gamma=b(b_\beta e^{-z_\gamma E_{-\gamma}}),\quad
b_\beta\in\tilde{a}_\beta{\text{N}}_\beta.
\end{equation*}
This is
\begin{equation*}
\begin{aligned}[t]
\,b_\gamma&=b(\tilde{a}_\beta e^{-z_\gamma
E_{-\gamma}})\in\tilde{a}_\gamma{\text{N}}_\gamma,\\
\tilde{a}_\gamma&:=\tilde{a}_\beta{\text{a}}^{H_\gamma}\,(\tilde{z}_\gamma\,%
)={\text{a}}^{H_\alpha}(\tilde{z}_\alpha){\text{a}}^{H_\beta}\,(\tilde{z}_%
\beta\,){\text{a}}^{H_\gamma}\,(\tilde{z}_\gamma\,),\\
\tilde{z}_\gamma&:=(\tilde{a}_\beta)^{-\gamma}z_\gamma={\text{a}}^{-\gamma(H_%
\alpha)/2}(\tilde{z}_\alpha){\text{a}}^{-\gamma(H_\beta)/2}\,(\tilde{z}_\beta%
\,){\text{a}}^{-\gamma(H_\beta)/2}\,(\tilde{z}_\gamma\,)z_\gamma
\end{aligned}
\end{equation*}
These are the desired formulas for these variables. The result is
\begin{equation*}
|c^\lambda_{{\text{s}}}(\xi)|^2=\prod{\text{a}}^{-\lambda_i(H_{\alpha_i})}(%
\tilde{z}_i)=\,a_\lambda(\xi)^{-1}
\end{equation*}
and
\begin{equation*}
\sigma_\lambda=-\frac{1}{2{\pi}{\text{i}}}\bar{\partial}\partial|c^\lambda_{{%
\text{s}}}|^2=+\frac{1}{2{\pi}{\text{i}}}\bar{\partial}\partial a_\lambda.
\end{equation*}
\end{proof}

\section{Canonical coordinates}

\noindent It will be convenient to change variables:
\begin{equation*}
\begin{aligned}[t]
&z_i=e^{\zeta_i}\quad\zeta_i:=\tau_i+\sqrt{-1}\phi_i\,\\
&\tilde{z}_i=e^{\tilde{\zeta}_i}\quad\tilde{\zeta}_i:=\tilde{\tau}_i+%
\sqrt{-1}\tilde{\phi}_i\,.
\end{aligned}
\end{equation*}
In these variables the formula for $\sigma_\lambda$ reads
\begin{equation*}
\begin{aligned}[t]
&\sigma_\lambda=\frac{1}{2{\pi}{\text{i}}}\sum_{ij}\frac{\partial^2K_%
\lambda}{\partial\bar{\zeta}_i\partial\zeta_j}d\bar{\zeta}_i\wedge
d\zeta_j=\frac{1}{2{\pi}}\sum_{ij}\frac{1}{2}\frac{\partial^2K_\lambda}{%
\partial\tau_i\partial\tau_j}d\tau_i\wedge d\phi_j,\\
&K_\lambda=\sum_il_i{\text{K}}(\tilde{\tau}_i)\quad{\text{K}}(\tilde{\tau}):=%
\log{\text{a}}(\tilde{\tau})\quad{\text{a}}(\tilde{\tau}):=1+e^{2\tilde{\tau}}
\end{aligned}
\end{equation*}
Also introduce
\begin{equation*}
J_i:=\frac{1}{2}\frac{\partial
K_\lambda}{\partial\tau_i},\quad{\text{J}}(\tilde{\tau}):=\frac{1}{2}\frac{d{%
\text{K}}(\tilde{\tau})}{d\tilde{\tau}}=\frac{e^{2\tilde{\tau}}}{1+e^{2%
\tilde{\tau}}}
\end{equation*}
The range of \text{J} is to be noted: $0\leq{\text{J}\leq 1.\text{ }}$ These
formulas have the following consequence.

\begin{theorem}
The curvature form of $\mathcal{L}_\lambda$ is given by 
\begin{equation*}
\sigma_\lambda=\frac{1}{2{\pi}}\sum_i\,dJ_i\wedge d\phi_i\,,
\end{equation*}
and the action of $e^H\in T_{{\text{s}}}$ is given by
\begin{equation*}
(J_i,\phi_i)\circ e^H=(J_i,\phi_i+\varpi_i(H)).
\end{equation*}
Define linear combinations $L_1,\cdots,L_\ell\,$of the $J_i$'s by 
\begin{equation*}
L_1:=0,\quad L_j:=\sum_{i:j>i}c_{ji}J_i.
\end{equation*}
Then $J_1,\cdots,J_\ell$ are recursively determined by
\begin{equation*}
J_1=l_1{\text{J}}_1,\quad
J_j={\text{J}}_j(l_j-L_j),\quad[{\text{J}}_j:={\text{J}}(\tilde{\tau}_j)]\quad
\end{equation*}
\ and their range is determined by $0\leq{\text{J}}_j\leq 1,\,$i.e.
\begin{equation*}
0\leq\frac{J_1}{l_1}\leq 1,\quad 0\leq\frac{J_j}{l_j-L_j}\leq 1
\end{equation*}
with the understanding that a vanishing denominator means a vanishing
numerator. 
\end{theorem}

\begin{proof}
The formulas for $\sigma_\lambda$ and the action of $T_{{\text{s}}}$ in terms
of the variables $(J_i,\phi_i)$ follows from the corresponding formulas in
terms of the variables $(z_i)$. The formula for $K_\lambda$ reads 
\begin{equation*}
\frac{1}{2}dK_\lambda=\sum_il_i{\text{J}}_id\tilde{\tau}_i,\quad{\text{J}}(%
\tilde{\tau}):=\frac{1}{2}\frac{d{\text{K}}(\tilde{\tau})}{d\tilde{\tau}}.
\end{equation*}
In terms of the $\tau_i$'s,
\begin{equation*}
\frac{1}{2}dK_\lambda=\sum_iJ_id\tau_i,\quad J_i:=\frac{1}{2}\frac{\partial
K_\lambda}{\partial\tau_i}\text{ .}
\end{equation*}
The relations between the $r_i=e^{\tau_i}$ and the
$\tilde{r}_i=e^{\tilde{\tau}_i}$ are
\begin{equation*}
r_\ell=\tilde{r}_\ell\quad\cdots\quad
r_i=\tilde{r}_i{\prod_{j:j>i}{\text{a}}}^{-\frac{1}{2}c_{ji}}(\tilde{r}_j).
\end{equation*}
In terms of the new variables these equations say
\begin{equation*}
\tau_i=\tilde{\tau}_i-\sum_{j:j>i}\frac{1}{2}c_{ji}{\text{K}}(\tilde{\tau}_j)
\end{equation*}
With some rearrangements this gives
\begin{equation*}
d\tau_i=d\tilde{\tau}_i+\sum_{j:j>i}c_{ji}{\text{J}}_jd\tilde{\tau}_j
\end{equation*}
with
\text{J}$_j:={\text{J}}(\tilde{\tau}_j):=(1/2)d{\text{K}}(\tilde{\tau})/$d$%
\tilde{\tau}$. Use these equations to compare the coefficients of the
$d\tilde{\tau}_j$'s in the two expressions 
\begin{equation*}
\frac{1}{2}dK_\lambda=\sum_jl_j{\text{J}}_jd\tilde{\tau}_j,\quad%
\frac{1}{2}dK_\lambda=\sum_iJ_id\tau_i\,.
\end{equation*}
Proceeding in the order $d\tilde{\tau}_1,d\tilde{\tau}_2,\cdots$ find
\begin{equation*}
l_1{\text{J}}_1=J_1,\quad
l_j{\text{J}}_j=J_j+\sum_{i:j>i}c_{ji}{\text{J}}_jJ_i
\end{equation*}
as required. 
\end{proof}

Some comments. (1)The theorem agrees with the results of Grossberg and
Karshon \lbrack 1994\rbrack, but their 2-form is constructed differently, by
a rather special procedure adapted to their setting. (2)If the form
$\sigma_\lambda$ is everywhere non-degenerate, the preceding theorem says
that the canonical coordinates $(J_i,\phi_i)$ are action-angle variables for
action of the compact real form of \ $T_{{\text{s}}}$ on the symplectic
manifold $\mathcal{X}_{{\text{s}}},\sigma_\lambda$ \lbrack Arnold,
1989\rbrack. \ In this context $\sigma_\lambda$ is the object of interest,
and indeed as form, not only as cohomology class. \ (3)In contrast to the
construction of the variables $\tilde{r}_\ell,\cdots,\tilde{r}_1$, the
inductive construction of $J_1,\cdots,J_\ell$ in terms of
$\tilde{\tau}_1,\cdots,\tilde{\tau}_\ell$ proceeds in the same direction as
the construction of the $\mathbb{P}^1$-chain. As a consequence,
$\,J_j={\text{J}}_j(l_j-\sum_{i:j>i}c_{ji}J_i)\,$depends only on the
variables $\tilde{\tau}_i$ with $i\leq j$ \ but on all of the variables
$\tau_i$.

For a curvature form of a holomorphic line bundle it is positivity which is
of primary interest, rather than just non-degeneracy. This is the question to
be addressed next.

\section{Positivity Questions}

\noindent Let \ $c_{ji}^\pm:=c_{ji}$ if $\pm c_{ji}>0$ and $=0\,$ otherwise.
\ Let $J_{j\text{ }\min},J_{j\text{ }\max}\,$ be the $\min/\max$ of $J_j$
$:={\text{J}}_j(l_j-\sum_{i:j>i}c_{ji}J_i)$ over $0\leq{\text{J}}_j\leq 1$
for fixed values of $J_{j-1},\cdots,J_1$, inductively given by \ 
\begin{equation*}
\begin{aligned}[t]
J_{j\min}&=\min\{0,\,l_j-\sum_{i:j>i}c_{ji}^{+}J_{i\text{
}\max}+c_{ji}^{-}J_{i\min}\}\\
J_{j\max}&=\max\{0,\,l_j-\sum_{i:j>i}c_{ji}^{+}J_{i\min}+c_{ji}^{-}J_{i\max}\}
\end{aligned}
\end{equation*}
\noindent\ \ $J_{j\text{ }\min\text{/}\max}=J_{j\text{
}\min\text{/}\max}(\lambda;$ $J_{j-1},\cdots,J_1$ $)$ are linear combinations
of $J_{j-1},\cdots,J_1$ depending on $\lambda$.

Let $\min J_j$ $,\max J_j\,$be the $\min/\max$ of $J_j$ over its complete
range, inductively given by
\begin{equation*}
\begin{aligned}[t]
\min J_j&=\min\{0,\,l_j-\sum_{i:j>i}c_{ji}^{+}\max J_i+c_{ji}^{-}\min J_i\}\\
\max J_j&=\max\{0,\,l_j-\sum_{i:j>i}c_{ji}^{+}\min J_i+c_{ji}^{-}\max J_i\}
\end{aligned}
\end{equation*}
\ $\min/\max J_j=\min/\max J_j(\lambda)$ depends on $\lambda$ only.

\begin{theorem}Assume $l_1\neq 0,\cdots,l_\ell\neq 0$. The form
$\sigma_{-\lambda}$ is positive if and only if $\lambda$ satisfies 
\begin{equation*}
\min J_1(\lambda)\geq 0\quad\cdots\quad\min J_\ell(\lambda)\geq 0.
\end{equation*}
If so, the range of $(J_i)$ is the cubic polytope given by the $\ell$
independent linear inequalities
\begin{equation*}
\begin{aligned}[t]
&\pi(\mathcal{X}):\quad 0\leq J_j\leq l_j-L_j,\quad
L_j:=\sum_{i:j>i}c_{ji}J_i.
\end{aligned}
\end{equation*}
\end{theorem}\noindent\textbf{Remarks.} \ (1) The conditions $\min
J_j(\lambda)\geq 0$ are equivalent to $\min J_j(\lambda)=0$ and amount to
$\ell$ piecewise linear inequalities in $\lambda$. They are certainly
satisfied if $l_\ell\gg l_{\ell-1}\gg\cdots\gg l_1>0;$ one only needs to
insure that so that the first term in $J_j$
$={\text{J}}_j(l_j-\sum_{i:j>i}c_{ji}J_i)$ dominates the rest.

(2)The action variables $J_j$ are determined only up to an additive constant
depending on $\lambda$ by the equation
$\sigma_\lambda=\frac{1}{2{\pi}}\sum_j\,dJ_i\wedge d\phi_j\,.\,$The $J_j$'s
used here are normalized so that $J_j|_{z_j=0}=0.\,$ Without this
normalization the positivity condition reads $\min
J_j(\lambda)=J_j|_{z_j=0}\,$. \ The polytope $\pi(\mathcal{X})$ is still
given by $J_{j\text{ }\min}\leq J_j\leq J_{j\text{ }\max},$ provided $0$ is
replaced by $J_j|_{z_j=0}$ in the definition of $J_{j\text{ }\min/\max}$.

The proof of the theorem needs some preparation. We have 
\begin{equation*}
J_j={\text{J}}_j(l_j-L_j),\quad L_j=\quad\sum_{i:j>i}c_{ji}J_i.
\end{equation*}
$L_j$ depends only on the variables$\,\tilde{\tau}_i$ with $j>i$ and
\begin{equation*}
{\text{J}}_i:={\text{J}}(\tilde{\tau}_i),\quad\text{
{K}}(\tilde{\tau}):=\log{\text{a}}(\tilde{\tau}),\quad{\text{J}}(\tilde{%
\tau}):=\frac{1}{2}\frac{d{\text{K}}(\tilde{\tau})}{d\tilde{\tau}}.
\end{equation*}
\begin{lemma}The curvature form of $\mathcal{L}_\lambda$ on
$\mathcal{X}_{{\text{s}}\text{ }}$is given by 
\begin{equation*}
\sigma_\lambda=\sum_{ji=1}^\ell\frac{1}{2}\frac{\partial^2K_\lambda}{\partial%
\tau_i\partial\tau_j}\frac{1}{4{\pi}{\text{i}}}\frac{d\bar{\zeta}_i\wedge
d\zeta_j}{\bar{\zeta}_i\zeta_j}\,.
\end{equation*}
Its $\ell$th exterior power is
\begin{equation*}
\sigma_\lambda^\ell=D_\lambda\,\prod_{j=1}^\ell\tilde{\sigma}_{-\alpha_j/2},%
\quad
D_\lambda:=\prod_j(l_j-L_j),\quad\tilde{\sigma}_{-\alpha_j/2}:=\prod_{j=1}^%
\ell\frac{{\text{K}}''(\tilde{\tau}_j)}{2{\pi}{\text{i}}}d\tilde{\tau}_j%
\wedge d\tilde{\phi}_j,
\end{equation*}
$\tilde{\sigma}_{-\alpha_j/2}$ being the curvature form on a copy of
$\mathbb{P}^1$ itself with affine coordinate
$\tilde{\zeta}_j=e^{\tilde{\tau}_j+\sqrt{-1}\tilde{\phi}_j}$.\end{lemma}

\begin{proof} (\emph{Lemma}) Write
$\sigma_\lambda=\frac{1}{2{\pi}{\text{i}}}\bar{\partial}\partial
K_\lambda\,$in terms of the operators
\begin{equation*}
\partial=\sum
d\zeta_j\,\frac{\partial}{\partial\zeta_j}\quad\bar{\partial}=\sum
d\bar{\zeta}_j\,\frac{\partial}{\partial\bar{\zeta}_j}\quad[d=\bar{\partial}+%
\partial].
\end{equation*}
Use $\,$the relations 
\begin{equation*}
dK_\lambda=2\sum J_jd\tau_j,\quad
d\tau_j=\frac{1}{2}(\partial\log\zeta_j+\bar{\partial}\log\bar{\zeta}_j)=%
\frac{1}{2}(\frac{d\zeta_j}{\zeta_j}+\frac{d\bar{\zeta}_j}{\bar{\zeta}_j})
\end{equation*}
to compute
\begin{equation*}
\sigma_\lambda=\frac{1}{2{\pi}{\text{i}}}\bar{\partial}\partial
K_\lambda=\frac{1}{{\pi}{\text{i}}}\sum\bar{\partial}J_j\wedge\partial\tau_j=%
\frac{1}{4{\pi}{\text{i}}}\sum\frac{\partial
J_j}{\partial\tau_i}\frac{d\bar{\zeta}_i\wedge
d\zeta_j}{\bar{\zeta}_i\zeta_j}\text{ .}
\end{equation*}
This proves the first assertion. To prove the second one, note that
\begin{equation*}
\det\frac{\partial J_j}{\partial\tau_i}=\det(\frac{\partial
J_j}{\partial\tilde{\tau}_l})\det(\frac{\partial\tilde{\tau}_l}{\partial%
\tau_i})=\prod_j\frac{\partial J_j}{\partial\tilde{\tau}_j}
\end{equation*}
since both determinants are triangular, the second with diagonal entries
$=1$. Since $L_j$ depends only on the variables$\,\tilde{\tau}_i$ with
$j>i\,,$
\begin{equation*}
\frac{\partial
J_j}{\partial\tilde{\tau}_j}=\frac{\partial((l_j-L_j){\text{J}}_j)}{\partial%
\tilde{\tau}_j}=(l_j-L_j){\text{J}}'(\tilde{\tau}_j)=\frac{1}{2}(l_j-L_j)\,{%
\text{K}}''(\tilde{\tau}_j).
\end{equation*}
Thus
\begin{equation*}
\det\frac{\partial
J_j}{\partial\tau_i}=\prod_j(l_j-L_j)\,\frac{1}{2}\,{\text{K}}''(\tilde{%
\tau}_j),\quad\sigma_\lambda^\ell=D_\lambda\,\prod_{j=1}^\ell\frac{{%
\text{K}}''(\tilde{\tau}_j)}{8{\pi}{\text{i}}}\frac{d\bar{\zeta}_j\wedge
d\zeta_j}{\bar{\zeta}_j\zeta_j}.
\end{equation*}
Since $\zeta_j=\tilde{\zeta}_j+$terms depending only on
$\tilde{\zeta}_{j+1},\cdots,\tilde{\zeta}_\ell,\,$the $\zeta_j$ can be
replaced by $\tilde{\zeta}_j$ in last formula. In terms of to
$\tilde{\tau}_j,\tilde{\phi}_j$ this gives 
\begin{equation*}
\sigma_\lambda^\ell=D_\lambda\,\prod_{j=1}^\ell\tilde{\sigma}_{-\alpha_j/2},%
\quad\tilde{\sigma}_{-\alpha_j/2}:=\prod_{j=1}^\ell\frac{{\text{K}}''(\tilde{%
\tau}_j)}{2{\pi}{\text{i}}}d\tilde{\tau}_j\wedge d\tilde{\phi}_j.
\end{equation*}
This proves the second assertion.\end{proof}

\begin{proof}(\emph{Theorem}) Generally, a (1,1) form
$\frac{1}{2{\text{i}}}\sum K_{ij}d\bar{z}_i\wedge dz_j$ is positive iff the
matrix $(K_{ij})$ is positive definite. So the form $\sigma_{-\lambda}$ is
positive definite at a point $\xi$ if and only if the matrix
$\partial^2K_{+\lambda}/\partial\tau_i\partial\tau_j$ is, and this \ happens
if and only if all principal minors of this matrix are positive \lbrack
Gantmacher, 1959, p.307\rbrack. \ Stated in terms of forms rather than
matrices, this criterion reads as follows. $\,$ \textit{A }$(1,1)$\textit{
form }$\sigma$\textit{ is positive at a point if and only if the sequence of
\ exterior powers }$\sigma^\ell,\sigma^{\ell-1},\cdots,\sigma$\textit{ is
positive on a complete flag in the tangent space}\emph{. }For a form of top
degree$,$ like $\sigma^\ell$, positivity means
$\sigma^\ell=D\prod\frac{1}{2{\text{i}}}d\bar{z}_j\wedge dz_j$ with $D>0,$ by
definition. 

To argue by induction on $\ell$, write the $\mathbb{P}^1$-chain
$\mathcal{X}_{{\text{s}}}$ as
$\mathcal{X}_{{\text{s}}}=G_\ell\times_{B_\ell}\mathcal{X}_{{\text{s}}'}\,%
$where $\mathcal{X}_{{\text{s}}'}$ is the $\mathbb{P}^1$-chain built from the
sequence \text{s}$'$ obtained from \text{s} by omitting the last letter
$i=\ell$. $\mathcal{X}_{{\text{s}}'}$ is embedded in
$\mathcal{X}_{{\text{s}}}$ as the fiber of
$\mathcal{X}_{{\text{s}}}\rightarrow G_\ell/B_\ell\approx\mathbb{P}^1$ over
the base-point 1$B_\ell$. The coordinate $z_\ell$ on
$\mathcal{X}_{{\text{s}}}$ is the pull-back of the affine coordinate on
centered at the point $s_\ell B$, also denoted $z_\ell$. The submanifold
$\mathcal{X}_{{\text{s}}'}$ of $\mathcal{X}_{{\text{s}}}$ is then given by
the equation $z_\ell\circ s_\ell=0$. Since $\tilde{\zeta}_\ell=\zeta_\ell$,
the formula for $\sigma_\lambda^\ell$ shows that at such a point
\begin{equation*}
\left.\sigma_\lambda^\ell\right|_{z_\ell\circ
s_\ell=0}=\,\,\left[(l_\ell-L_\ell)\frac{{\text{K}}''(\tau_\ell)}{8{\pi}{%
\text{i}}}\frac{d\bar{\zeta}_\ell\wedge
d\zeta_\ell}{\bar{\zeta}_\ell\zeta_\ell}\right]_{z_\ell\circ
s_\ell=0}\wedge\sigma_{\lambda'}^{\ell-1}
\end{equation*}
where $\lambda'$ is obtained by omitting the $\ell$th component and
$\sigma_{\lambda'}$ is the corresponding form on $\mathcal{X}_{{\text{s}}'}$.
The form in brackets is a pull-back of this form $G_\ell/B_\ell$ written in
the coordinate $\zeta_\ell$. It is positive if and only if $l_\ell-L_\ell$ is
positive. Since every point of $\mathcal{X}_{{\text{s}}}$ can be brought to
$\mathcal{X}_{{\text{s}}'}$ by the action of the compact form \text{U}$_\ell$
of \text{G}$_\ell$, which leaves $\sigma_\lambda$ invariant, one finds that
$\sigma_\lambda^\ell>0$ on $\mathcal{X}_{{\text{s}}}$ if and only if
\begin{equation*}
l_\ell-L_\ell>0\,\:{\text{and}}\:\,\sigma_{\lambda'}^{\ell-1}>0\,\:{%
\text{on}}\:\,\mathcal{X}_{{\text{s}}'}.
\end{equation*}
$l_\ell-L_\ell=l_\ell+\cdots$ is independent of $\tilde{\tau}_\ell$ and
$J_\ell={\text{J}}_\ell(l_\ell-L_\ell),\,0\leq{\text{J}}_\ell\leq 1$. Hence
$l_\ell-L_\ell>0$ amounts to $J_\ell\geq 0$ and not$\,\equiv 0$ on
$\mathcal{X}_{{\text{s}}}$. Since $J_\ell={\text{J}}_\ell l_\ell+\cdots$ the
latter condition is satisfied as soon as $\,l_\ell\neq 0,\,$which holds by
hypothesis, so that we may replace $l_\ell-L_\ell>0$ by $J_\ell\geq 0$. \
Thus $\sigma_{-\lambda}=-\sigma_\lambda$ is positive if and only if$\,J_1\geq
0,\,\cdots,J_\ell\geq 0$ throughout, i.e. $\min J_j(\lambda)=0$. If so, the
range of the$\,J_j={\text{J}}_j(l_j-L_j),0\leq{\text{J}}_j\leq 1$,is indeed
the given by $\,0\leq J_j\leq l_j-L_j$.\end{proof}

\noindent The condition for the positivity of $\sigma_\lambda$ will be
referred to as the \emph{positivity condition} on $\lambda$. \ 

\section{Quantization}

\noindent\textbf{Projective embedding}. Kodaira's theorems assert the
following consequences of positivity \lbrack Griffiths and Harris, 1978,
p.156, 181\rbrack.

\begin{theorem}Assume $\lambda$ satisfies the positivity condition.

\smallskip

\noindent\emph{(1)} Let $\Omega^q(\mathcal{L}_\lambda)$ be the sheaf of
holomorphic $q$-forms with values in $\mathcal{L}_\lambda$. 
\begin{equation*}
{\text{H}}^p(\mathcal{X}_{{\text{s}}},\Omega^q(\mathcal{L}_\lambda))=0\:{%
\text{for}}\:p+q>\ell.
\end{equation*}
\noindent\emph{(2)} Let $\mathcal{O}(\mathcal{L}_\lambda)$be the sheaf of
holomorphic sections of $\mathcal{L}_\lambda$. The natural map
\begin{equation*}
\mathcal{X}_{{\text{s}}}\rightarrow\mathbb{P}(\mathcal{H}^\ast),%
\mathcal{H}:={\text{H}}^0(\mathcal{X}_{{\text{s}}},\mathcal{O}(\mathcal{L}_{k%
\lambda})
\end{equation*}
is an embedding for any sufficiently large multiple $k\lambda\,$of $\lambda$.

\end{theorem}

Some comments. In general, the natural map
$\mathcal{X}\rightarrow\mathbb{P}(\mathcal{H}^\ast)$ into the projective
space dual to the space
$\mathcal{H}:={\text{H}}^0(\mathcal{X},\mathcal{O}(\mathcal{L}))$ of global
holomorphic sections is defined as follows. Given $\xi\in\mathcal{X}$, choose
any non-zero element $e(\xi)$ in the line $\mathcal{L}(\xi)$. \ The map
$f\mapsto f(\xi)/e(\xi)$ defines a linear functional
$e_\xi\in\mathcal{H}^\ast$, which up to a scalar depends only on $\xi$ and
represents the image of $\xi$ in $\mathbb{P}(\mathcal{H}^\ast)$. The line
bundle $\mathcal{L}$ on $\mathcal{X}$ is the pull-back of the hyperplane
bundle on the ambient projective space $\mathbb{P}(\mathcal{H}^\ast)$, so its
Chern class is the pull-back of that of the hyperplane bundle. These Chern
classes are represented by (normalized) curvature forms
$\sigma_\mathcal{X},\sigma_{\mathbb{P}^N}$ for any unitary metrics on these
line bundles. As form, however, $\sigma_\mathcal{X}$ need not coincide with
the pull-back of the Fubini-Study form $\sigma_{\mathbb{P}^N}$; \ on the
contrary, the relation between these forms is generally a delicate question.
\ In particular, the unitary norm on the Hilbert space of sections, which is
of interest for applications to representation theory, depends on the form
$\sigma_\mathcal{X}\,$itself. 

\medskip

\noindent\textbf{Some notation}. \ \ Throughout $\lambda$ is assumed to
satisfy the positivity condition and will remain fixed. It will be dropped as
a subscript, as will be the subscript \text{s}$,$ e.g. on $\sigma$ and
$\mathcal{X}$. Write $T$ for \ the image of $T_{{\text{s}}}$ as
transformation group on $\mathcal{X}$. Like any complex torus, $T$ has a
unique decomposition $T=T_{{\text{U}}}T_{{\text{P}}}$ into a 'unitary' real
form $T_{{\text{U}}}\approx\prod{\text{U}}(1)$ and a 'positive' real form
$T_{{\text{P}}}\approx{\prod\text{GL}(1,\mathbb{R})_{{\text{positive}}}}$. \
The weights $\varpi_i$ by which $T$ acts on the coordinates $z_i$ via
$z_i\circ h=h^{\varpi_i}z_i$ form a $\mathbb{Z}$--basis for a lattice
$\mathfrak{t}_{{\text{P}}}(\mathbb{Z})$ in the real dual \ of
$\mathfrak{t}_{{\text{P}}}^\ast$ of the of the Lie algebra of
$T_{{\text{P}}}$. Let $N_i$ be the dual basis for the dual lattice
$\mathfrak{t}_{{\text{P}}}(\mathbb{Z})$ in
$\mathfrak{t}_{{\text{P}}}(\mathcal{X})\subset\mathfrak{t}_{{\text{P}}}$,
defined by$\,\langle\varpi_i,N_j\rangle=\delta_{ij}$. \ The $N_i$ form an
$\mathbb{R}$-basis for $\mathfrak{t}_{{\text{P}}}$ and the \text{i}$N_i$ for
$\mathfrak{t}_{{\text{U}}}$. The imaginary unit \text{i}$:=\sqrt{-1}$ \ will
be written explicitly whenever $\mathfrak{t}_{{\text{U}}}$ is concerned: \
$\mathfrak{t}_{{\text{P}}}$ is taken as ``the'' real form of $\mathfrak{t}$
and $\mathfrak{t}=\mathfrak{t}_{{\text{P}}}+\mathfrak{t}_{{\text{U}}}$ as the
decomposition of $\mathfrak{t}$ into ``real'' and ``imaginary'' parts. For
example, write \ 
\begin{equation*}
{\text{i}}H\in\mathfrak{t}_{{\text{U}}}:H=\sum_i\epsilon_iN_i,\epsilon_i\in%
\mathbb{R},\quad{{\text{i}}}\eta\in\mathfrak{t}_{{\text{U}}}^\ast:\eta=%
\sum_in_i\varpi_i\,,\,n_i\in\mathbb{R}.
\end{equation*}
\ \ Thus \text{i} functions as map
$\mathfrak{t}_{{\text{P}}}\rightarrow\mathfrak{t}_{{\text{U}}}$. \ The
pairing of $H\in\mathfrak{t}_{{\text{P}}}$ and
$\eta\in\mathfrak{t}_{{\text{P}}}^\ast$ is written as $\langle\eta,H\rangle$
or $\langle H,\eta\rangle$. \medskip

\noindent\textbf{The moment set}. Any element
\text{i}$H={\text{i}}\sum\epsilon_iN_i$ of $\mathfrak{t}_{{\text{U}}}$ has
$J_{{\text{i}}H}(\xi):=2{\pi}\sum\epsilon_iJ_i(\xi)$ as a Hamiltonian
function, in the sense that $dJ_{{\text{i}}H}=-\iota({\text{i}}H)\sigma$, the
inner product of $\sigma$ with the vector field on $\mathcal{X}$ induced by \
$-{\text{i}}H\in\mathfrak{t}_{{\text{U}}}$. This follows from the formulas
$\sigma=\frac{1}{2{\pi}}\sum dJ_i\wedge
d\phi_i,\,\iota({\text{i}}N_i)dJ_j=0,$
$\iota({\text{i}}N_i)d\phi_j=\delta_{ij}.\,$In particular the Hamiltonian
function for \text{i}$N_i$ is the action variable $J_i(\xi)$. Since
$J_{{\text{i}}H}(\xi)$ depends linearly on \text{i}$H$, it defines a map
$\mathcal{X}\rightarrow\mathfrak{t}_{{\text{U}}}^\ast,$ the \emph{moment map}
of the $T_{\text{U}}$ action on $\mathcal{X}$. \ Transferred to
$\mathfrak{t}_{{\text{P}}}^\ast$ \ it becomes a map
$\pi:\mathcal{X}\rightarrow\mathfrak{t}_{{\text{P}}}^\ast$, defined by
$\langle\pi(\xi),H\rangle:=J_{{\text{i}}H}(\xi)$ i.e. $\pi(\xi)=\sum
J_i(\xi)\varpi_i$. \ The image $\pi(\mathcal{X})$ of $\pi$ in
$\mathfrak{t}_{{\text{P}}}^\ast$ can evidently be identified with the range
of the $J_i$ and will be called the \emph{moment set}. \medskip

\noindent\textbf{The weight set.} \ Let $\mathcal{H}$ be the space of
holomorphic sections of $\mathcal{L},$ a finite-dimensional space with a
positive definite norm defined by the formula
\begin{equation*}
\|f\|^2:=\int_\mathcal{X}|f|^2\,\sigma^\ell.
\end{equation*}
There is a natural representation of the complex torus on $\mathcal{H}$,
defined by 
\begin{equation*}
(U_hf)(\xi):=f(h^{-1}\xi),\,
\end{equation*}
which is unitary on the real torus $T_{{\text{U}}}$. \ 

Any \text{i}$\eta\in\mathfrak{t}_{{\text{U}}}^\ast(\mathbb{Z})$ defines a
Laurent monomial $z^\eta:=z_1^{n_1}\cdots z_\ell^{n_\ell}$ in the coordinates
$z_i=e^{\tau_i+{\text{i}}\phi_i}$ on $\mathcal{X}$ and hence a holomorphic
section $f_\eta=z^\eta c^\lambda$ on the coordinate domain. ($f_\eta$
transforms by $\tilde{\eta}:=\eta+\tilde{\lambda}$ under $T_{{\text{U}}}$
since \ $c^\lambda$ transforms by $\tilde{\lambda}$.) \ The \emph{weight set}
$\Pi\,$ is the set of weights
$\eta\in{\text{i}}^{-1}\mathfrak{t}_{{\text{U}}}^\ast(\mathbb{Z})$ for which
$f_\eta$ extends to a holomorphic section on all of $\mathcal{X}$, a
condition which may be reformulated in other ways.

\begin{lemma}Let $f=z^\eta c^\lambda$. \ The following conditions are
equivalent.

\emph{(1)} $f(\xi)$ extends to a global holomorphic section on \
$\mathcal{X}$.

\emph{(2)}$|f(\xi)|$ is bounded on $\mathcal{X}$.

\emph{(3)} $|f(\xi)|^2\sigma^\ell$ is integrable over $\mathcal{X}$.

\end{lemma}

\begin{proof} It suffices to prove the equivalence locally, i.e. with
$\mathcal{X}$ replaced by a neighbourhood of each of its points. Since the
coordinates are meromorphic functions on $\mathcal{X}$, $f$ is a meromorphic
section on $\mathcal{X}$, represented locally by a meromorphic function,
regular of on coordinate domain, and for such a function the local conditions
are indeed equivalent.\end{proof}

\begin{theorem} The weight set $\Pi$ $\subset\mathfrak{t}_{{\text{P}}}^\ast$
is $\Pi=$
$\pi(\mathcal{X})\cap{\text{i}}^{-1}\mathfrak{t}_{{\text{U}}}^\ast(%
\mathbb{Z})$, i.e. $f_\eta=z^\eta c^\lambda$ extends to a global holomorphic
section on all of $\mathcal{X}$ if an$d$ only if \ $\eta$ belongs to the
image $\pi(\mathcal{X})$ of the moment map. \ These
$f_\eta,\,\eta\in\Pi,\,$form a basis for $\mathcal{H}$.\end{theorem}

\begin{proof} The coordinates $z_i=e^{\tau_i+{\text{i}}\phi_i}$ on
$\mathcal{X}$ provide a bijection 
\begin{equation*}
T_{{\text{U}}}\backslash\mathcal{X}\leftrightarrow\mathfrak{t}_{\text{P}},%
\quad T_{{\text{U}}}\xi\mapsto H:=\sum\tau_i(\xi)N_i,
\end{equation*}
defined on the domain where all coordinates are non-zero. Equivalently,
$\xi=e^H\xi_1$ where $\xi_1$ has coordinates $z_1=\cdots=z_\ell=1$. \ The
pointwise square norm $a(\xi)=|c^\lambda(\xi)|^2$ is really a function on
$T_{{\text{U}}}\backslash\mathcal{X}$. Use the bijection
$T_{{\text{U}}}\xi\mapsto H$ to write the equation $K(\xi)=\log a(\xi)$ as
$K(H)=\log a(H)$ in terms of $H\in\mathfrak{t}_{{\text{P}}}$. 

By the lemma, $f_\eta$ extends iff its pointwise square-norm is bounded. This
square norm is 
\begin{equation*}
|f_\eta|(\xi)^2=|z^\eta
c^\lambda(\xi)|^2=|z|^{2\eta}e^{-K(\xi)}=e^{2\langle\eta,H\rangle-K(H)}\text{
}
\end{equation*}
where $H=\sum_i\epsilon_iN_i,\,\epsilon_i=\tau_i(\xi)$. Thus the condition is
precisely that $2\langle\eta,H\rangle-K(H)$ be bounded from above as function
of $H$. \ Since the curvature form 
\begin{equation*}
\sigma=\sum_{ji=1}^\ell\frac{1}{2}\frac{\partial^2K(H)}{\partial\tau_i%
\partial\tau_j}\frac{1}{4{\pi}{\text{i}}}\frac{d\bar{\zeta}_i\wedge
d\zeta_j}{\bar{\zeta}_i\zeta_j}\,.
\end{equation*}
is positive definite, so is the Hessian matrix of $K(H)$. Hence $K(H)$
\textit{is strictly convex} as function of $H\in\mathfrak{t}_{{\text{P}}}$
\lbrack Rockafellar, 1970, p.27\rbrack. \ For a smooth convex \ function
$K(H)$, $2\langle\eta,H\rangle-K(H)$ achieves a finite maximum iff its
partials vanish at some point $H_o$ \lbrack loc. cit., \ p.258\rbrack, and
this amounts \ to \ $2n_i=\partial K/\partial\tau_i=2J_i$ at $H_o$, i.e.
$\eta\in\pi(\mathcal{X})$. \ The elements $f_\eta,\eta\in\Pi,\,$ are
precisely the weight vectors for the representation of $T_{{\text{U}}}$ on
$\mathcal{H},\,$hence form a basis. \end{proof}

\label{a}\begin{corollary}The function
$Z({\text{i}}H):={\text{trace}}(e^{{\text{i}}H}\mid\mathcal{H})$ of
$e^{{\text{i}}H}\in T_{{\text{U}}}\,$is given by the formula
\begin{equation*}
Z({\text{i}}H)=\sum_{\eta\in\Pi}e^{{\text{i}}\eta(H)}\,.
\end{equation*}
\end{corollary}\noindent A comment. The determination of the weight set by
convex calculus, as in the proof of the theorem, is as general as it is
simple. Useful tools are available, especially the apparatus of \ conjugate
functions. \ By definition, the \textit{conjugate} of a convex function
$K(H)$ is $K^{\text{c}}(\eta):=\max\{2\langle\eta,H\rangle-K(H)\mid
H\in\mathfrak{t}_{{\text{P}}}\}$ as function on the dual space
$\mathfrak{t}_{{\text{P}}}^\ast$. (In this context $\max=\infty$ is allowed,
by convention.) In the language of convex analysis, the weight set could be
described as the set of integral points in the domain of the conjugate
function $K^{\text{c}}(\eta),\eta\in\mathfrak{t}_{\text{P}}^\ast$ of the
K\"ahler potential $K(H),H\in\mathfrak{t}_{{\text{P}}}$. To an analyst convex
conjugates will be familiar as asymptotic phase functions according to the
principle of stationary phase, \lbrack Guillemin and Sternberg, 1977,
p.397\rbrack.\medskip

\noindent\textbf{Partition functions}. Whether the function $Z({\text{i}}H)$
is to be called \textit{(quantum) partition function} or \emph{character} is
a matter of taste and tradition. The prejudice displayed here is adopted for
this occasion only, because of a curious relation between toric geometry and
statistical mechanics, as will now be explained .

Fix \text{i}$H\in\mathfrak{t}_{{\text{U}}}$ and let $J:=J_{{\text{i}}H}$ be
its Hamiltonian function. The classical density associated to this
Hamiltonian is the function on $\mathcal{X}$ by the formula$\,$
\begin{equation*}
\rho_{\text{classical}}(\xi):=\frac{1}{Z_{\text{classical}}}\frac{\,e^{-\beta
J(\xi)}}{\ell!},\quad\int_\mathcal{X}\rho_{\text{classical}}\,\sigma^\ell=1
\end{equation*}
\lbrack Honerkamp 1998, p.28\rbrack. \ The denominator $Z_{\text{classical}}$
is a normalizing factor, depending only on $H$: 
\begin{equation*}
Z_{\text{classical}}:=\int_\mathcal{X}e^{-\beta J}\frac{\sigma^\ell}{\ell!}.
\end{equation*}
It can be written in the form 
\begin{equation*}
Z_{\text{classical}}:=\int_\mathcal{X}e^{-\beta J+\sigma}
\end{equation*}
the exponential being taken in the exterior algebra with the understanding
that the integral picks out the component of the relevant degree
$2\ell=\dim_\mathbb{R}\mathcal{X}$. 

On the other hand, the quantum density is the density in the spectral
resolution of \ $e^{-\beta{\text{i}}H}$ as operator on $\mathcal{H}$:$\,$
\begin{equation*}
\rho_{\text{quantum}}(\eta):=\frac{m(\eta)}{Z_{\text{quantum}}}\,e^{-\beta{%
\text{i}}\eta},\quad\sum_\eta\rho_{\text{quantum}}(\eta)=1,
\end{equation*}
$m(\eta)$ being the multiplicity of \ the \text{i}$\eta$ as eigenvalue of
\text{i}$H$ on $\mathcal{H}$. The normalizing factor $Z_{\text{quantum \
}}$is therefore \ 
\begin{equation*}
Z_{\text{quantum }}=\sum_\eta
m(\eta)e^{-\beta{\text{i}}\eta}={\text{Tr}}(e^{-\beta{\text{i}}H}\mid%
\mathcal{H}).
\end{equation*}
\lbrack Honerkamp, 1998, p.204\rbrack. To bring out the dependence on
\text{i$H\in\mathfrak{t}_{{\text{U}}}$} write $\rho({\text{i}}H;\eta)$ and
$Z({\text{i}}H)$, if necessary.

The fundamental relation between the classical and quantum partition function
is \textbf{Kirillov's Formula} \lbrack Berline, Getzler, and Vergne,
p.250\rbrack. In the present situation it states that for \text{i}$H$ in a
neighbourhood of zero in $\mathfrak{t}_{{\text{U}}},$ on which the
exponential map is one-to-one, one has 
\begin{equation*}
\int_\mathcal{X}e^{J_{{\text{i}}H}+\sigma}\widehat{A}({\text{i}}H)=\sum_\eta
m(\eta)e^{{\text{i}}\langle\eta,H\rangle}.
\end{equation*}
The sum on the right is the quantum partition function and the integral on
the left is the classical partition function, except for the extra factor
$\widehat{A}({\text{i}}H)$ under the integral, a characteristic class
involving the Riemann curvature of $\mathcal{X}$, whose definition can be
found loc. cit. \ (Its appearance here suggests that the classical definition
of the partition function should be modified so as to incorporate this
factor, a modification which is irrelevant in a flat phase space.) That
\text{i}$H$ has to remain in neighbourhood of zero is a serious restriction.

Returning to the toric setting, assume that Kodaira's map
$\psi:\mathcal{X}\rightarrow\mathbb{P}(\mathcal{H}^\ast),\xi\mapsto\psi_\xi,%
$is an embedding. \ Let $\mathcal{X}'$ be the image of $\mathcal{X}$ and $T'$
the image of $T$ as transformation group of $\mathcal{X}'$. Let
$\mathcal{L}'$ be the line bundle on $\mathcal{X}'$ induced by the hyperplane
bundle on $\mathcal{H}^\ast$. Let $\sigma'$ be the curvature from for the
unitary metric on $\mathcal{L}'$ induced by any $T_{\text{U}}$-invariant
unitary metric on $\mathcal{H}^\ast$, for example the metric defined above.
Let $J'_{{\text{i}}H'}(\xi)$ be a Hamiltonian function for
\text{i}$H'\in\mathfrak{t}'_{{\text{U}}}$ acting on $(\mathcal{X}',\sigma')$,
i.e. \ $dJ'_{{\text{i}}H'}=-\iota({\text{i}}H')\sigma'$. This equation
determines $J'_{{\text{i}}H'}(\xi)$ up to a constant depending linearly on
\text{i}$H'$. \ Fix this constant as follows. Let $\xi_1$ be a base-point for
$\mathcal{X}$ in the open $T$-orbit, say the point \ with all coordinates
$z_i$ equal to 1. Let $\xi'_1$ be it image in $\mathcal{X}'$, and require
that $J'_{{\text{i}}H'}(\xi'_1)=J_{{\text{i}}H}(\xi_1)$ if
\text{i}$H'\in\mathfrak{t}_{{\text{U}}}'$ is the image of
\text{i}$H\in\mathfrak{t}_{{\text{U}}}$. Thus we have another set of data
$\mathcal{X}',\sigma'\cdots\,$like $\mathcal{X},\sigma\cdots$. 

\begin{lemma}Under the identification $\mathfrak{t}'=\mathfrak{t}$, 
\begin{equation*}
Z_{\text{classical}}'({\text{i}}H)=Z_{\text{classical}}({\text{i}}H)\text{,}%
\quad Z_{\text{quantum}}'({\text{i}}H)=Z_{\text{quantum}}({\text{i}}H).
\end{equation*}
\end{lemma}

\begin{proof}The equality of the quantum partition functions is obvious,
since the line bundles $\mathcal{L}$ and $\mathcal{L}'$ correspond under
$\mathcal{X}=\mathcal{X}'.\,$ The equality of the classical partition
functions, however, is not, since the forms $\sigma$ and $\sigma'$ need not
coincide. What is obvious is that their cohomology classes coincide, since
both represent the Chern class of $\mathcal{L}'=\mathcal{L}$. \
Thus$\,\sigma'=\sigma+d\varphi$ for some smooth 1-form $\varphi$, which may
be chosen $T_{{\text{U}}}$-invariant like $\sigma$ and $\sigma'$. \ Then
$d\circ\iota({\text{i}}H)\varphi=-\iota({\text{i}}H)\circ d\varphi$ for all
\text{i}$H\in\mathfrak{t}_{{\text{U}}}$ and (suppressing \text{i}$H$) the
equations \ $dJ'+\iota\sigma'=0,\,dJ+\iota\sigma=0,\,$give \
$dJ'=dJ+d\iota\varphi,\,$i,e $J'=J+\iota\varphi+C$. The normalization
$J'(\xi_1)=J(\xi_1)$ implies $C=0$, so
$J'+\sigma'=\,J+\sigma+(d+\iota)\varphi$. \ The rest of the proof uses three
observations: \ (1)$(d+\iota)(J+\sigma)=0$, (2)$(d+\iota)$ is a derivation,
and (3)$(d+\iota)\Phi=d\Phi+$forms of lower degree. Armed with these one
computes 
\begin{equation*}
\begin{aligned}[t]
\int_\mathcal{X}e^{J+\sigma}e^{(d+\iota)\varphi}&=\int_\mathcal{X}(e^{J+%
\sigma}+e^{J+\sigma}(d+\iota)\Phi))\quad[e^{(d+\iota)\varphi}=:1+(d+\iota)%
\Phi]\\
&=\int_\mathcal{X}(e^{J+\sigma}+(d+\iota)(e^{J+\sigma}\Phi))\quad[\,{{%
\text{by}\:\text{ (1)}\:\text{ and}\:\text{ (2)}}}\,]\\
&=\int_\mathcal{X}e^{J+\sigma}+0+0\quad[\,{\text{by}\:\text{ Stokes}\:\text{
and}\:\text{ (3)}}\,]
\end{aligned}
\end{equation*}
This equation says $Z'_{{\text{classical}}}=Z_{{\text{classical}}}$ except
for the factor $\beta$, which may be absorbed into $J$ and $J'$. \end{proof}

The image of the measure $|\sigma^\ell|=|\prod\frac{1}{2{\pi}}dJ_i\wedge
d\phi_i|$ \ under $\eta=\pi(\xi)$ i.e.$\,n_i=J_i(\xi)$ is evidently
$\rho(\eta)d\eta$, the Haar measure $d\eta=$ $\prod dn_i$ multiplied by the
indicator function $\rho(\eta)$concentrated on $\pi(\mathcal{X})$. \ After an
integration over the angle variables the partition function becomes 
\begin{equation*}
Z_{\text{classical}}({\text{i}}H)=\int_{\eta\in\pi(\mathcal{X})}e^{-\beta%
\langle\eta,H\rangle}d\eta,
\end{equation*}
the Laplace transform of $\rho(\eta)$. \ The lemma therefore implies that the
moment sets for $\mathcal{X}'$ and $\mathcal{X}$ coincide, i.e. \
$\pi'(\mathcal{X})=\pi(\mathcal{X})$ under the identification
$\mathfrak{t}'=\mathfrak{t}$. 

The main point is the topological nature of classical partition function:
$Z_{{\text{classical}}}(H)$ \emph{depends only on the Chern class of the line
bundle.} (This \ phenomenon is evident in the setting of Kirillov's formula
in the form presented loc. cit, where the lemma is understood in terms of
equivariant cohomology. The calculation of Fourier transforms like one above
was one of the first applications of that theory \lbrack Berline and Vergne,
1983\rbrack.) Replacing $\sigma$ by $\sigma'$ allows one to apply the theory
of toric varieties, which provides information to the moment map and the
weight set. There is in particular an intriguing description of the moment
set in terms of the maximum entropy principle of statistical mechanics, due
to C. Lee \lbrack 1990, p.13, Ewald 1996, p.298\rbrack.

\medskip

\noindent\textbf{The path integral formula. \ }The unitary representation of
$T_{{\text{U}}}$ on $\mathcal{H}$ is a quantization of the symplectic
manifold $\mathcal{X},\sigma$ associated to the totally complex, positive
polarization provided by its complex structure \lbrack Woodhouse,
1991\rbrack; the wave functions are the holomorphic sections $f(\xi)$ which
make up $\mathcal{H}$. The quantum evolution operators are those representing
$T_{{\text{U}}}$. On the other hand, the canonical variables
$(J,\phi):=(J_i,\phi_i)$ provide a natural real polarization, at least on
their domain of definition, whose wave functions are functions of $\phi$. (It
natural to choose $\phi$ rather than $J$, because $(J,\phi)$ has an
interpretation as coordinates on the cotangent bundle of $T$ with $\phi$ as
coordinate on the base.) We now leave the realm of geometric quantization for
that of path integral quantization. The evolution operator generated by a
Hamiltonian function $H(J,\phi)$ is then supposed to be an integral operator
on wave functions $\psi(\phi)$,
\begin{equation*}
U_H\psi(\phi'')=\int K_H(\phi'',\phi')\psi(\phi')d\phi',
\end{equation*}
whose kernel is given by a path integral, formally written as 
\begin{equation*}
K_H(\phi'',\phi'):=\int e^{\frac{{\text{i}}}{\hbar}S_H[J(t),\phi(t)]}.
\end{equation*}
$S_H[J(t),\phi(t)]\,$is the action along a path $(J(t),\phi(t)\mid t'\leq
t\leq t'')$, defined by \ 
\begin{equation*}
S_H[J(t),\phi(t)]:=\int_{t'}^{t''}[Jd\phi-H(J,\phi)dt].
\end{equation*}
The path integral for $K_H(\phi'',\phi')$ is supposed to be taken over all
such path with $\phi(t')=\phi',\phi(t'')=\phi''$. The quantum partition
function $Z({\text{i}}H)$, now written without subscript, is the trace of
$U_H$, obtained from $K_H(\phi'',\phi')$ by setting $\phi''=\phi'\,$and
integrating. The difficulties encountered when trying to define path
integrals as true integrals are well known and need not be discussed here.
Instead, I shall give a recipe for evaluating $Z({\text{i}}H)$, which can be
taken as a definition of this path integral in the special case we are
concerned with, i.e. for a $\mathbb{P}^1$-chain $\mathcal{X}$, equipped with
its curvature form $\sigma$ and canonical coordinates $(J,\phi)$, and for
Hamiltonians which are linear functions $H=H(J)$ of the $J$'s only.

Consider first the case of $\mathbb{P}^1$ itself. The coordinates are
$(J,\phi)\,$with $J$ ranging over an interval $\pi(\mathcal{X}):J_{\min}\leq
J\leq J_{\max}\,$ and $\phi$ over $\mathbb{R}/2{\pi}\mathbb{Z}$. The
Hamiltonian is $H(J):=HJ$ for some constant also denoted $H$. \emph{Define}
the path integral
\begin{equation*}
Z({\text{i}}H)=\int_{\phi(0)\equiv\phi(1)}e^{{\text{i}}S_H[J(t),\phi(t)]}
\end{equation*}
by the following recipe.

\medskip

\noindent(1) Fix $N$ and replace $S_H[J(t),\phi(t)]$ by the sum 
\begin{equation*}
\begin{aligned}[t]
&S^{(N)}:=[J\phi]_0^1-\sum_{i=0}^{N-1}\phi_i\Delta J_i+HJ_i\Delta t\text{,}\\
&\Delta J_i:=J_{i+1}-J_i,\quad\Delta t:=1/N.
\end{aligned}
\end{equation*}
\noindent(2) Fix $n\in\mathbb{Z}$ and integrate over all sequences
$(J_0\phi_0,\cdots,J_N\phi_N)$
\begin{equation*}
J_{\min}\leq J_i\leq
J_{\max},\quad-\infty<\phi_i<\infty,\quad\phi_N-\phi_0=2{\pi}n.
\end{equation*}

\noindent(3) Sum over $n\in\mathbb{Z}$.

\medskip

\noindent(4) Take the limit $N\rightarrow\infty$.

\medskip

\noindent Written as a single formula (not recommended), this amounts to the
definition
\begin{equation*}
\begin{aligned}[t]
&Z({\text{i}}H):=\lim_{N\rightarrow\infty}\sum_{n=0}^\infty\int_{J_{%
\min}}^{J_{\max}}dJ_N\int_{-\infty}^\infty
d\phi_N\,\cdots\,\,\int_{J_{\min}}^{J_{\max}}dJ_1\int_{-\infty}^\infty
d\phi_1\times\\
&\times\delta(\phi_N-\phi_0-2{\pi}n)e^{i(J_N\phi_N-J_0\phi_0)-i%
\sum_{i=1}^{N-1}\phi_i(J_i-J_{i-1})+HJ_i(\frac{1}{N})}.
\end{aligned}
\end{equation*}
The integrals in this formula must be interpreted in the sense of generalized
functions, i.e. as limits of integrals against compactly supported test
functions tending to the constant function $\equiv 1$. The integrals over
$J_{\max}\leq J\leq J_{\min},$ in particular, are defined as limits against
continuous compactly supported functions identically $\equiv 1$ on this
closed interval. In (1), an integration by parts has been applied to the
integral for $S_H$, replacing the form $Jd\phi$ by $-\phi dJ$ before
replacing $J(t),\phi(t)$ by $(J_1\phi_1,\cdots,J_N\phi_N)$. In (2), the
variables $\phi_i$ replacing $\phi(t)\,$are allowed to vary from $-\infty$ to
$+\infty$, while the endpoint condition
$\phi(0)\equiv\phi(1)\,{\text{mod}\,}2{\pi}$ is interpreted as the \emph{sum}
over all such paths and \emph{not} by requiring that the angle itself return
to its initial value. This is an essential feature.

In the case of a general toric $\mathbb{P}^1$-chain the coordinate pair
$(J,\phi)$ is replaced by $\ell$ coordinate pairs$\,(J_i,\phi_i)$ and the
interval $\pi(\mathcal{X})=\{J_{\max}\leq J\leq J_{\min}$ $\}$ by the cubic
polytope $\pi(\mathcal{X})=\{J_{i\text{ }\min}\leq J_i\leq J_{i\text{
}\max}\}$. The path integral $Z({\text{i}}H)$ is then defined as an
$\ell$-fold integral, each of which is given by the recipe above. 

\begin{theorem}The path integral for $Z({\text{i}}H)$ exists and is given by
\begin{equation*}
Z({\text{i}}H)=\sum_{\eta\in\Pi}\,e^{{\text{i}}\langle\eta,H\rangle}\,.
\end{equation*}
\end{theorem}

\begin{proof}Follow the steps in the definition: 

\medskip

\noindent(1 ) Fix $N$. As indicated let,
\begin{equation*}
\begin{aligned}[t]
S_H&=\int_0^1(Jd\phi-HJdt)=[J\phi]_0^1-\int_0^1(\phi dJ+HJdt)\\
S_H^{(N)}&:=[J\phi]_0^1-\sum_{k=0}^{N-1}[\phi_k(J_{k+1}-J_k)+\frac{1}{N}HJ_k]
\end{aligned}
\end{equation*}
\noindent(2) Fix $n$. For fixed $k$, the part of the integral involving
$(J_k\phi_k)$ is of the form
\begin{equation*}
\int_{J_{\min}}^{J_{\max}}dJ_k\int_{-\infty}^\infty
d\phi_k\,e^{-{\text{i}}\phi_k(J_{k+1}-J_k)-\frac{{\text{i}}}{N}HJ_k\,}
\end{equation*}
Use the Fourier inversion formula on the real line in the form 

\begin{equation*}
\begin{aligned}[t]
&\int_{-\infty}^\infty\frac{d\phi_k}{2{\pi}}\,\,e^{-{\text{i}}%
\phi_k(J_{k+1}-J_k)}=\delta(J_{k+1}-J_k),\\
&\int_{-\infty}^\infty dJ_k\,f(J_k)\delta(J_k-J_{k-1})=f(J_{k-1}).
\end{aligned}
\end{equation*}
This gives
\begin{equation*}
\begin{aligned}[t]
Z({\text{i}}H)^{(N,n)}&:=\int_{J_{\min}}^{J_{\max}}dJ_{N\,}\int_{J_{%
\min}}^{J_{\max}}dJ_{0\,}\,e^{i(J_N\phi_N-J_0\phi_0)-iHJ_0}\,\delta(J_N-J_0)\\
&=\int_{J_{\min}}^{J_{\max}}d\eta\,e^{-{\text{i}}\eta(\phi_N-\phi_0+H)}
\end{aligned}
\end{equation*}

\noindent(3)By Poisson summation, 
\begin{equation*}
\sum_{n\in\mathbb{Z}}Z({\text{i}}H)^{(N,n)}=\sum_{n\in\mathbb{Z}}\int_{J_{%
\min}}^{J_{\max}}d\eta\,e^{-{\text{i}}\eta(2{\pi}n+H)}=\sum_{\eta\in\Pi}\,%
e^{{\text{i}}\eta H}.
\end{equation*}
(4)The result of (3) is independent of $N$; the limit $N\rightarrow\infty$ is
trivial.

\medskip

\noindent The proof has been formulated so that it remains valid without
change in the setting of $\mathbb{P}^1$-chains. \end{proof}

Some comments. (1) The appearance of the delta function $\delta(J_N-J_0)$ in
step (2) may be interpreted to mean that the path integral for
$Z({\text{i}}H)$ is automatically concentrated on paths $(J(t),\phi(t))$
closing in $J$ as well as in $\phi$, which is not part of the input.

(2) The theorem gives $Z({\text{i}}H)$ again as a sum over the set of
integral points in the range of the action variables $(J_i)$. Since the proof
is nothing but Poisson's summation formula applied to the definition of the
path integral, as specified in the recipe, the canonical coordinates
themselves remain as the essential ingredient. \ The applicability of the
path integral in many other situations seems to indicate that it does capture
some basic feature, and not only in a formal way.

(3) In a more philosophical vein, it appears that the unitary-positive
splitting $T=T_{{\text{U}}}T_{{\text{P}}}$ is reflected in the statistical
mechanics vs. quantum mechanics setting of the two formulas for the partition
functions. The odd nature of this juxtaposition is an old puzzle \lbrack
Feynman and Hibbs, 1965, p.296\rbrack.

\end{document}